\documentclass[english,times]{elsarticle}

\usepackage[margin=2cm]{geometry}
\usepackage{lineno,hyperref}
\usepackage{color}
\usepackage{float}
\usepackage{bm}
\usepackage{amsmath}
\usepackage{amssymb}
\usepackage{amsthm}
\usepackage{graphicx}
\usepackage{setspace}
\usepackage{subscript}
\usepackage{upgreek}
\usepackage{natbib}
\usepackage{babel}
\usepackage[algo2e,ruled]{algorithm2e} 

\newcommand{\diag}{\mathop{\mathrm{diag}}}

\DeclareMathOperator{\tr}{tr}

\newtheorem*{remark}{Remark}
\newtheorem{theorem}{Theorem}
\newtheorem{corollary}{Corollary}
\SetKwInOut{Input}{Input}
\SetKwInOut{Output}{Output}
\SetKwInOut{Initialize}{Initialize}
\DeclareMathAlphabet{\mathpzc}{OT1}{pzc}{m}{it}

\newcommand{\RN}[1]{%
  \textup{\uppercase\expandafter{\romannumeral#1}}%
}

\journal{Computers $\&$ Chemical Engineering}

\begin{document}
\begin{frontmatter}

\title{Stochastic data-driven model predictive control using Gaussian processes}

\author[mymainaddress]{Eric Bradford \corref{mycorrespondingauthor1}}

\author[mymainaddress]{Lars Imsland}

\author[mymbinaddress,mymcinaddress]{Dongda Zhang}

\author[mymbinaddress]{Ehecatl Antonio del Rio Chanona
\corref{mycorrespondingauthor2}}

\cortext[mycorrespondingauthor1]{Corresponding author: eric.bradford@ntnu.no}

\cortext[mycorrespondingauthor2]{Corresponding author: a.del-rio-chanona@imperial.ac.uk}

\address[mymainaddress]{Department of Engineering Cybernetics, Norwegian University of Science and Technology, Trondheim, Norway}

\address[mymbinaddress]{Centre for Process Systems Engineering (CPSE), Department of Chemical Engineering, Imperial College London, UK}

\address[mymcinaddress]{School of Chemical Engineering and Analytical Science, University of Manchester, Manchester, UK}

\cortext[]{\\ \textcopyright 2020 The Authors. This manuscript version is made available under the CC-BY-NC-ND 4.0 license http://creativecommons.org/licenses/by-nc-nd/4.0/. This is a post-peer-review, pre-copyedit version of an article published in the journal Computers $\&$ Chemical Engineering (CACE). The final authenticated version is available online at: https://doi.org/10.1016/j.compchemeng.2020.106844.}

\begin{abstract}
Nonlinear model predictive control (NMPC) is one of the few control methods that can handle multivariable nonlinear control systems with constraints. Gaussian processes (GPs) present a powerful tool to identify the required plant model and quantify the residual uncertainty of the plant-model mismatch. It is crucial to consider this uncertainty, since it may lead to worse control performance and constraint violations. In this paper we propose a new method to design a GP-based NMPC algorithm for finite horizon control problems. The method generates Monte Carlo samples of the GP offline for constraint tightening using back-offs. The tightened constraints then guarantee the satisfaction of chance constraints online. Advantages of our proposed approach over existing methods include fast online evaluation, consideration of closed-loop behaviour, and the possibility to alleviate conservativeness by considering both online learning and state dependency of the uncertainty. The algorithm is verified on a challenging semi-batch bioprocess case study. 
\end{abstract}

\begin{keyword}
Model-based nonlinear control, Uncertain dynamic systems, Machine learning, Probabilistic constraints, State space, Robust control
\end{keyword}

\end{frontmatter}

\section{Introduction}
Model predictive control (MPC) describes an advanced control method that has found a wide range of applications in industry. MPC employs an explicit dynamic model of the plant to determine a finite sequence of control actions to take at each sampling time. The main advantage of MPC is its ability to deal with multivariate plants and process constraints explicitly \citep{Maciejowski2002}. Linear MPC is relatively mature and well-established in practice, however many systems display strong nonlinear behaviour motivating the use of nonlinear MPC (NMPC) \citep{Allgower2004}. NMPC is becoming progressively more popular due to the advancement of improved non-convex optimization algorithms \citep{Biegler2010}, in particular in chemical engineering \citep{Biegler1991}.

The performance of MPC is however greatly influenced by the accuracy of the plant model, which is one of the main reasons why MPC is not more widely used in industry \citep{Lucia2014}. The development of an accurate plant model has been cited to take up to 80\% of the MPC commissioning effort \citep{Sun2013}. NMPC algorithms exploit various types of models, commonly developed by first principles or based on process mechanisms \citep{Nagy2007b}. Many mechanistic and empirical models are however often too complex to be used online and in addition have often high development costs. Alternatively, black-box identification models can be exploited instead, such as support vector machines \citep{Xi2007}, fuzzy models \citep{Kavsek-Biasizzo1997}, neural networks (NNs) \citep{Piche2000}, or Gaussian processes (GPs) \citep{Kocijan2004}. For example, recently in \citet{Wua,Wu} recurrent NNs are utilised for an extensive NMPC approach with proofs on closed-loop state boundedness and convergence applied to a chemical reactor. In addition, in \citet{Wu2019a} the approach was extended updating the recurrent NNs online to further improve the effectiveness.            

GPs are an interpolation technique developed by \citet{Krige1951} that were popularized by the machine learning community \citep{Rasmussen2006}. While GPs have been predominantly used to model static nonlinearities, there are several works that apply GPs to model dynamic systems \citep{Girard2003,Kocijan2005,Bradford2018e}. GP predictions are given by a Gaussian distribution. The mean of this distribution can be viewed as a deterministic prediction, while the variance can be interpreted as a measure of uncertainty for this deterministic prediction. This uncertainty measure is generally difficult to obtain by nonlinear parametric models \citep{Kocijan2004} and may in part explain the relative popularity of GPs. For control applications this uncertainty measure can be utilized to efficiently learn a dynamic model by exploring unknown regions or avoiding regions with too high uncertainty to improve robustness \citep{Berkenkamp2015}. So far, GPs have been exploited in a multitude of ways in the control community, including reinforcement learning \citep{Deisenroth2011}, designing robust linear controllers \citep{Umlauft2017}, or adaptive control \citep{Chowdhary2015}. In particular, GPs have been shown to be an efficient approach to attain approximate plant models for NMPC.                        

The use of GPs for NMPC was first proposed in \citet{Murray-Smith2003}, which updates a GP model online for reference tracking without constraints. In \citet{Kocijan2004,Kocijan2005a} the GP is instead identified offline and utilized online, in which the variance is constrained to prevent the NMPC from steering the dynamic system into regions with high uncertainty. A GP plant model is updated online in \citet{Klenske2016} and in \citet{Maciejowski2013} to overcome unmodeled periodic errors or changes to the dynamic system after a fault has occurred respectively. GPs have been shown in \citet{Wang2016} to be an efficient means for disturbance forecasting for a linear stochastic MPC approach applied to a drinking water network. GPs have further been applied to approximate the mean and variance required in stochastic NMPC \citep{Bradford2018b}. \citet{Grancharova2007} derived an explicit solution for GP-based NMPC. In \citet{Cao2017} a GP dynamic model is employed for the control of an unmanned quadrotor, while in \citet{Likar2007} a GP dynamic model is exploited to control a gas-liquid separation process. While these and other works show the feasibility of GP-based MPC, there is a lack of efficient approaches to account for the uncertainty measure provided. Model uncertainty can lead to constraint violations and worse performance. To mitigate the effect of uncertainty on MPC, robust MPC \citep{Bemporad1999a} and stochastic MPC \citep{Mesbah2016,Heirung2018} methods have been developed.      

The majority of works for GP-based MPC indeed consider the uncertainty measure provided, however most proposed algorithms employ stochastic uncertainty propagation to achieve this, for example \citep{Kocijan2004,Kocijan2005a,Hewing2018,Cao2017,Grancharova2007,Wang2016}. \citet{Hewing2017} give an overview of the various stochastic propagation techniques available. These approaches have some considerable disadvantages, which are:
\begin{itemize} \itemsep0em  
    \item{No known methods to exactly propagate stochastic uncertainties through GP models. Instead, only approximations are available relying on linearization or statistical moment-matching.}
    \item{Increased computational time of GP-based MPC due to the propagation approach itself.}
    \item{Most works consider only open-loop propagation of uncertainties, which is often prohibitively conservative due to open-loop growth of uncertainties.}
\end{itemize}

Recently some papers have proposed different robust GP-based MPC algorithms. In \citet{Koller2018a} a NMPC algorithm is introduced based on propagating ellipsoidal sets using linearization, that provides closed-loop stability guarantees. This approach may however suffer from increased computational times, since the ellipsoidal sets are propagated online. Furthermore, the method may be relatively conservative due to the use of Lipschitz constants. \citet{Maiworm2018} propose the use of a robust MPC approach by bounding the one-step ahead error, while the determination of the required parameters seems to be relatively difficult. \citet{Soloperto2018} suggest a robust control approach for linear systems, in which the GP is used to represent unmodeled nonlinearities. The approach is shown to stabilize the linear system despite these uncertainties, which however may have no solution if the difference between the linearized system and the actual nonlinear system is too large.     

In this paper we extend an algorithm first introduced in \citet{Bradford2019b}, for which the following extensions were made:
\begin{itemize} \itemsep0em
    \item Inclusion of uncertainty for the initial state.
    \item Adding additive disturbance noise to the problem definition.
    \item Accounting for state dependency on the GP noise.
    \item Improved algorithm to obtain the required back-offs using root-finding as opposed to the inverse CDF, which has superior convergence and leads to improved satisfaction of the required probability bounds. 
\end{itemize}

The aim of this approach is to take into account the uncertainty given by a GP state space model for a NMPC \textit{finite-horizon} control problem. Due to the issues using stochastic uncertainty propagation for the NMPC formulation as highlight previously, we base the NMPC only on cheap evaluations of the GP. This leads to considerably faster evaluation times with little effect on the performance. The proposed method utilizes explicit back-offs, which were recently proposed in \cite{Koller2018,Paulson2018} to account for stochastic uncertainties in NMPC. These methods generally rely on generating closed-loop Monte Carlo (MC) samples offline from the plant to attain the required back-off values. To obtain exact MC samples of the GP dynamic models we exploit results from \citet{Conti2009,Umlauft2018}. There are several important advantages of this new method: 

\begin{itemize} \itemsep0em
    \item{Back-offs are attained using closed-loop simulations, therefore the issue of open-loop growth of uncertainties is avoided.}
    \item{Required computations are carried out offline, such that the online computational times are nearly unaffected.}
    \item{Independence of samples allows some probabilistic guarantees to be given.}
    \item{Explicit consideration of online learning and state dependency of the uncertainty to alleviate conservativeness.}
\end{itemize}

The proposed method is a data-driven NMPC approach relying on black-box identification of a GP model from input/output data pairs. The required data may be obtained either by simulations from a high-fidelity model or by experiments using for example step tests. In the algorithm the state dependency of the uncertainty is accounted for by introducing a penalty term on the variance in the objective. This variance for GPs is a function of the states and hence leads to a trade-off between exploiting the GP model to optimize the objective and avoiding uncertain regions to reduce the spread of the trajectories. The algorithm proposed is aimed at \textit{finite horizon} control problems, for which batch processes are a particularly important example. They are utilized in many different chemical engineering sectors due to their inherent flexibility to deal with variations in feedstock, product specifications, and market demand. Frequent highly nonlinear behaviour and unsteady-state operation of batch processes have led to the increased acceptance for advanced control solutions, such as NMPC \citep{Nagy2007a}. Works on batch process NMPC accounting for uncertainties include an extended and Unscented Kalman filter based algorithms for uncertainty propagation \citep{Nagy2003,Bradford2018c}, a NMPC algorithm using min-max successive linearization \citep{Valappil2002}, NMPC algorithms that employ PCEs to account for possible parametric uncertainties \citep{Mesbah2014,Bradford2019d}, and multi-stage NMPC \citep{Lucia2013}.      

The paper is comprised of the following sections. In Section \ref{sec:prob_def} the problem definition is given. In Section \ref{sec:Gaussian_processes} a general outline of GPs is given including the sampling procedure used. Section \ref{sec:Solution_approach} shows how the GPs can be exploited to solve the defined problem. In Section \ref{sec:case_study} the semi-batch bioprocess case study is described, while in Section \ref{sec:results_discussions} the results and discussions for this case study are given. Section \ref{sec:conclusions} concludes the paper.    

\subsection*{Notation}
$\mathbb{N}$ and $\mathbb{R}$ represent the sets of natural numbers and real numbers respectively. The variable $\delta_{ij}$ denotes the Kronecker delta function, such that:
\begin{equation*}
    \delta_{ij} :=
\begin{cases}
    1 , & \text{if} \, \, i = j\\
    0 , & \text{otherwise}
\end{cases}
\end{equation*}

The notation $\diag(a_0,a_1,\ldots,a_n)$ is used to represent the following diagonal matrix:
\begin{equation*}
    \diag(a_0,a_1,\ldots,a_n) := \begin{bmatrix} 
a_0    & 0      & \ldots & 0        \\
0      & a_1    & \ddots & \vdots   \\
\vdots & \ddots & \ddots & 0        \\
0      & \ldots & 0      & a_n   
\end{bmatrix}
\end{equation*}

We represent the Gaussian distribution with mean $\bm{\upmu}$ and covariance $\bm{\Sigma}$ as $\mathcal{N}(\bm{\upmu},\bm{\Sigma})$. Further, $\bm{\upphi} \sim \mathcal{N}(\bm{\upmu},\bm{\Sigma})$ denotes that the random variable $\bm{\upphi}$ follows a Gaussian distribution with mean $\bm{\upmu}$ and covariance $\bm{\Sigma}$. 

The expected value of a random variable $\bm{\upphi}$ is denoted as:
\begin{equation*}
    \mathbb{E}[\bm{\upphi}] := \int_{\Omega} \bm{\upphi} dp_{\bm{\upphi}} 
\end{equation*}
where $p_{\bm{\upphi}}$ the probability density function of $\bm{\upphi}$ over the sample space $\Omega$.

Further, we define the indicator function and the probability measure of random variable $\bm{\upphi}$ as follows:
\begin{align*}
&     \mathbf{1}\{C \leq c\} := \begin{cases}
    1 , & \text{if} \, \, C \leq c \\
    0 , & \text{otherwise}
\end{cases} \\
& \mathbb{P}\{\bm{\upphi} \in \mathcal{A} \} := \int_{\bm{\upphi} \in \mathcal{A}} \bm{\upphi} dp_{\bm{\upphi}}, \quad \mathbb{P}\{\bm{\upphi} \leq c\} := \mathbb{E}[\mathbf{1}\{\bm{\upphi} \leq c\}] 
\end{align*}
where $\mathcal{A}$ is a set defining an event on $\bm{\upphi}$ and $\mathbb{P}\{\bm{\upphi} \leq c\}$ is the probability that $\bm{\upphi}$ is less than or equal to $c$.  

Lastly, we require the definition of the beta inverse cumulative distribution function (cdf) for a random variable $\phi$. This function $\text{betainv}(P,A,B)$ returns a value $C$ of $\upphi$ following a beta distribution with parameters $A,B$ that has a probability of $P$ to be less than or equal to $C$. The definition is as follows:
\begin{align*}
    & \text{betainv}(P,A,B) \in F_{\phi}^{-1}(P|A,B) = \{ C:F_{\phi}(C|A,B)=P) \} \\
    & F_{\phi}(C|A,B) := \frac{1}{\mathcal{B}(A,B)} \int_{0}^{C} t^{A-1} (1-t)^{B-1} dt, \quad \mathcal{B}(A,B) = \int_{0}^1 t^{A-1} (1-t)^{B-1} dt  
\end{align*}

\section{Problem definition} \label{sec:prob_def}
The dynamic system in this paper is given by a discrete-time nonlinear equation system with additive disturbance noise:
\begin{equation}
    \mathbf{x}_{t+1} = \mathbf{f}(\mathbf{x}_t,\mathbf{u}_t) + \bm{\upomega}_t, \quad \mathbf{x}_0 \sim \mathcal{N}(\bm{\upmu}_{\mathbf{x}_0},\bm{\Sigma}_{\mathbf{x}_0}) \label{eq:f}
\end{equation}
where $t$ is the discrete time, $\mathbf{x} \in \mathbb{R}^{n_{\mathbf{x}}}$ is the state, $\mathbf{u} \in \mathbb{R}^{n_{\mathbf{u}}}$ are the control inputs, $\mathbf{f}:\mathbb{R}^{n_{\mathbf{x}}} \times \mathbb{R}^{n_{\mathbf{u}}} \rightarrow \mathbb{R}^{n_{\mathbf{x}}}$ are nonlinear equations, and $\bm{\upomega}$ represents Gaussian distributed additive disturbance noise with zero mean and diagonal covariance matrix $\bm{\Sigma}_{\bm{\upomega}}$. The initial condition $\mathbf{x}_0$ is assumed to be Gaussian distributed with mean $\bm{\upmu}_{\mathbf{x}_0}$ and covariance matrix $\bm{\Sigma}_{\mathbf{x}_0}$. 

We assume measurements of the states to be available with additive Gaussian noise expressed as:
\begin{equation}
    \mathbf{y} = \mathbf{f}(\mathbf{x},\mathbf{u}) + \bm{\upnu} \label{eq:y}
\end{equation}
where $\mathbf{y} \in \mathbb{R}^{n_{\mathbf{x}}}$ is the measurement of $\mathbf{f}(\mathbf{x},\mathbf{u})$ perturbed by additive Gaussian noise $\bm{\upnu} \sim \mathcal{N}(\mathbf{0},\bm{\Sigma}_{\bm{\upnu}})$ with zero mean and a diagonal covariance matrix $\bm{\Sigma}_{\bm{\upnu}}=\diag(\sigma^2_{\nu_1},\ldots,\sigma^2_{\nu_{n_{\mathbf{x}}}})$.  

The aim of the control problem is to minimize a finite-horizon cost function:
\begin{equation}
    V_T(\mathbf{x}_0,\mathbf{U}) = \mathbb{E} \left[\sum_{t=0}^{T-1} \ell(\mathbf{x}_t,\mathbf{u}_t) + \ell_f(\mathbf{x}_T) \right] \label{eq:objdef}
\end{equation}
where $T \in \mathbb{N}$ is the time horizon, $\mathbf{U}=[\mathbf{u}_{0},\ldots,\mathbf{u}_{T-1}]^{\sf T} \in \mathbb{R}^{ T \times n_{\mathbf{u}}}$ is a joint vector of $T$ control inputs, $\ell:\mathbb{R}^{n_{\mathbf{x}}} \times \mathbb{R}^{n_{\mathbf{u}}} \rightarrow \mathbb{R}$ is the stage cost, and $\ell_f:\mathbb{R}^{n_{\mathbf{x}}} \rightarrow \mathbb{R}$ denotes the terminal cost.   

The control problem is subject to hard constraints on the inputs:
\begin{align}
    & \mathbf{u}_t \in \mathbb{U}_t && \forall t \in \{0,\ldots,T-1\} \label{eq:ucon}
\end{align}

The states are subject to a joint chance constraint that requires the satisfaction of a nonlinear constraint set up to a certain probability, which can be stated as: 
\begin{subequations} \label{eq:xcon}
\begin{align}
    & \mathbb{P} \left\{ \bigcap^T_{t=0} \{ \mathbf{x}_t \in \mathbb{X}_t \}    \right\} \geq 1-\epsilon
\end{align}
where $\mathbb{X}_t$ is defined as:
\begin{align}
    & \mathbb{X}_t = \{ \mathbf{x} \in \mathbb{R}^{n_{\mathbf{x}}} \mid g_j^{(t)}(\mathbf{x}) \leq 0, j=1,\ldots,n_g \}
\end{align}
\end{subequations}
The joint chance constraints are formulated such that the joint event over all $t \in \{0,\ldots,T\}$ of all $\mathbf{x}_t$ fulfilling the nonlinear constraint sets $\mathbb{X}_t$ has a probability greater than $1-\epsilon$.

For convenience we define the tuple $\mathbf{z} = (\mathbf{x},\mathbf{u}) \in \mathbb{R}^{n_{\mathbf{z}}}$ with joint dimension $n_{\mathbf{z}}=n_{\mathbf{x}} + n_{\mathbf{u}}$. The dynamic system in Equation \ref{eq:f} is assumed to be unknown. Instead, we are only given a \textit{finite number} of noisy measurements according to Equation \ref{eq:y}. The available data can then be denoted by the following two matrices:
\begin{subequations} \label{eq:datasetdef}
\begin{align}
    & \mathbf{Z} = [\mathbf{z}^{(1)},\ldots,\mathbf{z}^{(N)}]^{\sf T} \in \mathbb{R}^{N \times n_{\mathbf{z}}} \\
    & \mathbf{Y} = [\mathbf{y}^{(1)},\ldots,\mathbf{y}^{(N)}]^{\sf T} \in \mathbb{R}^{N \times n_{\mathbf{x}}} 
\end{align}
\end{subequations}
where $\mathbf{z}^{(i)}$ represents the input of the $i$-th data point with corresponding noisy observation $\mathbf{y}^{(i)}$, $N$ denotes the overall number of training data points, $\mathbf{Z}$ is a collection of input data, and the corresponding noisy observations are collected in $\mathbf{Y}$.

It should be noted that the uncertainty in this problem arises partially from the uncertain initial condition $\mathbf{x}_0$ and the additive disturbance noise $\bm{\upomega}$. Most of the uncertainty however comes from the fact that we do not know $\mathbf{f}(\mathbf{x},\mathbf{u})$ and are only given noisy observations of $\mathbf{f}(\mathbf{x},\mathbf{u})$ instead. To solve this problem we train a GP to approximate $\mathbf{f}(\cdot)$ using the available data in Equation \ref{eq:datasetdef}. The GP methodology is introduced for this purpose in the next section. This GP then represents a distribution over possible functions $\mathbf{f}(\cdot)$ given the available data, which can be exploited to attain stochastic constraint satisfaction of the closed-loop system.

\section{Gaussian processes} \label{sec:Gaussian_processes}

\subsection{Regression} \label{sec:Gaussian_processes_regression}
In this section we introduce the use of GPs to infer a latent function $f:\mathbb{R}^{n_{\mathbf{z}}} \rightarrow \mathbb{R}$ from noisy data. For a more complete overview refer to \citet{Rasmussen2006}. Let the noisy observations $y$ of $f(\cdot)$ be given by:
\begin{align}
    & y = f(\mathbf{z}) + \nu \label{eq:yGP}
\end{align}
where $\mathbf{z} \in \mathbb{R}^{n_{\mathbf{z}}}$ is the argument of $f(\cdot)$ and $y$ is a perturbed observation of $f(\mathbf{z})$ with additive Gaussian noise $\nu \sim \mathcal{N}(0,\sigma^2_{\nu})$ with zero mean and variance $\sigma^2_{\nu}$.

GPs can be considered a generalization of multivariate Gaussian distributions to describe a distribution over functions. A GP is fully specified by a mean function and a covariance function. The mean function represents the "average" shape of the function, while the covariance function specifies the covariance between any two function values. We write that a function $f(\cdot)$ is distributed as a GP with mean function $m(\cdot)$ and covariance function $k(\cdot,\cdot)$ as:
\begin{align}
    & f(\cdot) \sim \mathcal{GP}(m(\cdot),k(\cdot,\cdot)) \label{eq:fGPdist}
\end{align}

The prior GP distribution is defined by the \textit{choice} of the mean function and covariance function. In this study we apply a zero mean function and the squared-exponential (SE) covariance function defined as:  
\begin{align}
    & m(\mathbf{z}) \coloneqq 0 \label{eq:meanfundef}   \\
    & k(\mathbf{z},\mathbf{z}') \coloneqq \zeta^2 \exp \left( -\frac{1}{2}(\mathbf{z} - \mathbf{z}')^{\sf T}\bm{\Lambda}^{-2}(\mathbf{z} - \mathbf{z}') \right) \label{eq:covfundef}
\end{align}
where $\mathbf{z},\mathbf{z}' \in \mathbb{R}^{n_{\mathbf{z}}}$ are arbitrary inputs, $\zeta^2$ denotes the covariance magnitude, and $\bm{\Lambda}^{-2} \coloneqq \diag(\lambda_1^{-2},\ldots,\lambda_{n_{\mathbf{z}}}^{-2})$ is a scaling matrix. 

\begin{remark}[Prior assumptions]
Note the zero mean assumption can be easily achieved by normalizing the data beforehand. The SE covariance function is smooth and stationary, such that choosing the SE covariance function assumes the latent function $f(\cdot)$ to be smooth and stationary as well. The algorithm presented in this work can be utilised using any covariance function. In the case of highly non-stationary functions it may be necessary to use non-stationary covariance functions \citep{Sampson1992}.     
\end{remark}

From the additive property of Gaussian distributions the measurements of $f(\cdot)$ also follow a GP accounting for measurement noise:
\begin{align}
    & y \sim \mathcal{GP}(m(\mathbf{z}),k(\mathbf{z},\mathbf{z}') + \sigma^2_{\nu} \delta_{\mathbf{z} \mathbf{z}'}) \label{eq:fGPdisty}
\end{align}

We denote the hyperparameters defining the prior jointly by $\bm{\Psi} \coloneqq [\zeta,\lambda_1,\ldots,\lambda_{n_{\mathbf{z}}},\sigma_{\nu}]^{\sf T}$, in which the variance $\sigma_{\nu}$ of the measurement noise is included in case it is unknown. Commonly the hyperparameters are unknown, such that they need to be inferred from the available data using for example maximum likelihood estimation (MLE). 

Assume we are given $N$ noisy function evaluations according to Equation \ref{eq:yGP} denoted by $\mathbf{Y} \coloneqq [y^{(1)},\ldots,y^{(N)}]^{\sf T} \in \mathbb{R}^N$ as the result of the inputs given in $\mathbf{Z}=[\mathbf{z}^{(1)},\ldots,\mathbf{z}^{(N)}]^{\sf T} \in \mathbb{R}^{N \times n_{\mathbf{z}}}$. According to the prior GP assumption, the data follows a multivariate Gaussian distribution:
\begin{align}
    & \mathbf{Y} \sim \mathcal{N}(\mathbf{0},\bm{\Sigma}_{\mathbf{Y}})
\end{align}
where $[\bm{\Sigma}_{\mathbf{Y}}]_{ij}=k(\mathbf{z}^{(i)},\mathbf{z}^{(j)}) + \sigma^2_{\nu} \delta_{ij}$ for each pair $(i,j) \in \{1,\ldots,N\}^2$.

The log-likelihood of the observations is consequently given by (ignoring constant terms):
\begin{align}
    & \mathcal{L}(\bm{\Psi}) \coloneqq -\frac{1}{2} \log(\det(\bm{\Sigma}_{\mathbf{Y}})) - \frac{1}{2} \mathbf{Y}^{\sf T} \bm{\Sigma}_{\mathbf{Y}}^{-1} \mathbf{Y} \label{eq:MLE}
\end{align}

The MLE estimate of the hyperparameters $\bm{\Psi}$ is determined by maximizing Equation \ref{eq:MLE}. Once the hyperparameters are known, we need to determine the posterior GP distribution of the latent function $f(\cdot)$. From the prior GP assumption we know that the training data and the value of $f(\cdot)$ at an arbitrary input $\mathbf{z}$ follow a joint multivariate normal distribution:
\begin{align}
    & \begin{bmatrix}
    \mathbf{Y}   \\
    f(\mathbf{z}) 
\end{bmatrix} \sim \mathcal{N}\left(\begin{bmatrix}
    \mathbf{0}   \\
    0 
\end{bmatrix},\begin{bmatrix}
    \bm{\Sigma}_{\mathbf{Y}} & \mathbf{k}^{\sf T}(\mathbf{z})  \\
    \mathbf{k}(\mathbf{z}) & k(\mathbf{z},\mathbf{z}')
\end{bmatrix}\right) \label{eq:jointnormal}
\end{align}
where $\mathbf{k}(\mathbf{z}) \coloneqq [k(\mathbf{z},\mathbf{z}^{(1)}),\ldots,k(\mathbf{z},\mathbf{z}^{(N)})]^{\sf T}$.

The posterior Gaussian distribution of $f(\mathbf{z})$ given the data $(\mathbf{Z},\mathbf{Y})$ can then be found by using the conditional distribution rule for multivariate normal distributions based on the joint normal distribution in Equation \ref{eq:jointnormal}, which leads to: 
\begin{subequations}
\begin{align}
    & f(\mathbf{z}) | \mathcal{D} \sim \mathcal{N}(\mu_f(\mathbf{z};\mathcal{D}),\sigma_f(\mathbf{z};\mathcal{D})) \label{eq:postdist}
\end{align}
with 
\begin{align}
& \mu_f(\mathbf{z};\mathcal{D})    \coloneqq \mathbf{k}^{\sf T}(\mathbf{z}) \bm{\Sigma}_{\mathbf{Y}}^{-1} \mathbf{Y} \label{eq:postmean} \\
& \sigma_f^2(\mathbf{z};\mathcal{D}) \coloneqq \zeta^2 - \mathbf{k}^{\sf T}(\mathbf{z}) \bm{\Sigma}_{\mathbf{Y}}^{-1} \mathbf{k}(\mathbf{z}) \label{eq:postvar}
\end{align}
\end{subequations}
where $\mathcal{D}=(\mathbf{Z},\mathbf{Y})$ denotes the training data available to obtain the posterior Gaussian distribution. The mean function $\mu_f(\mathbf{z};\mathcal{D})$ in this context is the prediction of the GP at $\mathbf{z}$, while the variance function $\sigma_f^2(\mathbf{z};\mathcal{D})$ is a measure of uncertainty.   

\begin{figure}[H] \centering
   \includegraphics[width=0.8\textwidth]{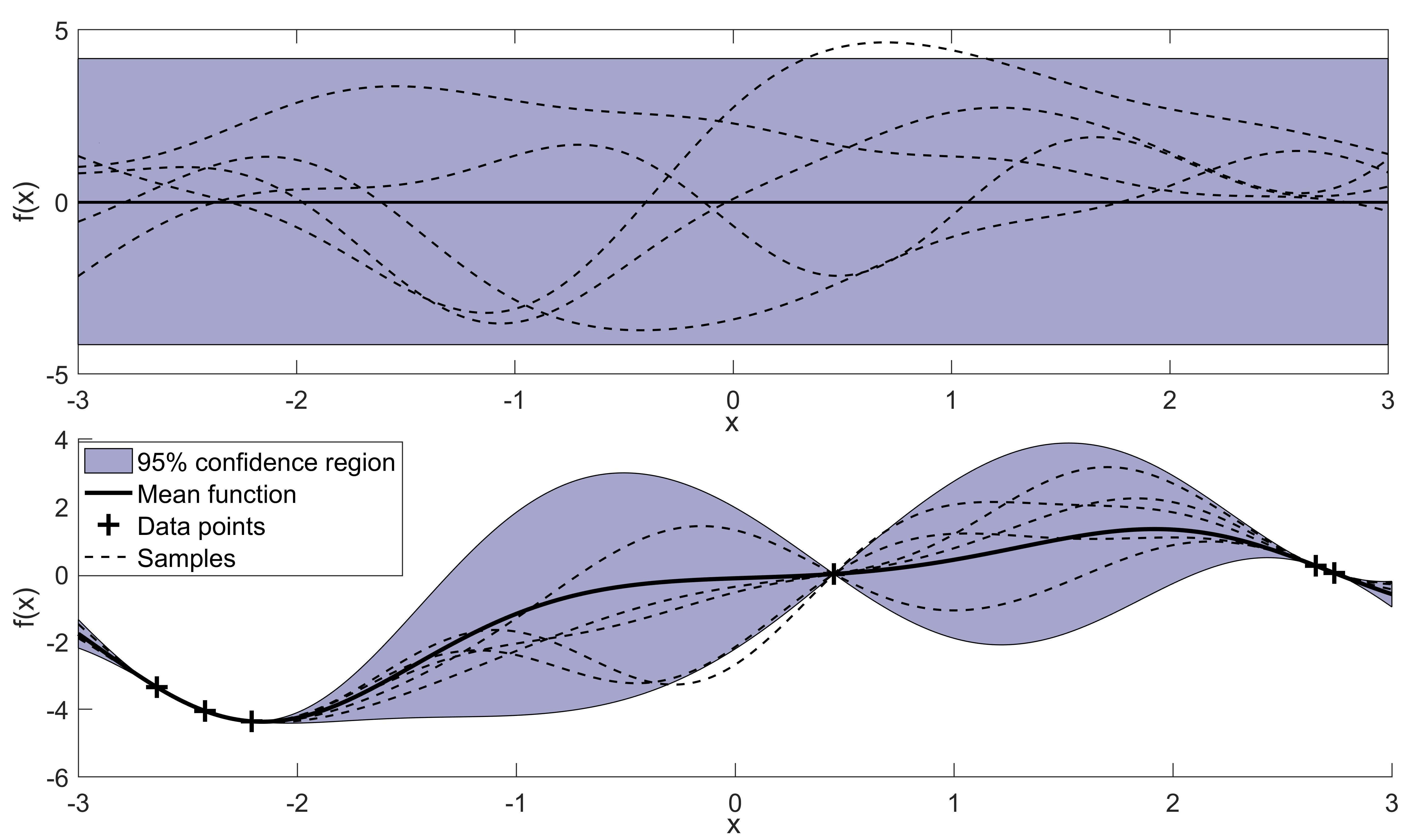}
  \caption{Illustration of a GP of a 1-dimensional function perturbed by noise. On the top the prior of the GP is shown, while on the bottom the Gaussian process was fitted to several observations to obtain the posterior.}
  \label{fig:Prior_posterior}
\end{figure}

In Figure \ref{fig:Prior_posterior} we illustrate a prior GP in the top graph and the posterior GP in the bottom graph.

\subsection{Recursive update} \label{sec:Gaussian_processes_online_learning}
Often we are given a training dataset $\mathcal{D}$ initially to build a posterior GP and then obtain data points individually afterwards. For example in GP-based MPC we may build an initial GP offline following the procedure shown in Section \ref{sec:Gaussian_processes_regression}, and then also update this model online as new data becomes available. In this work we keep the hyperparameters constant but update the mean function and variance function in Equations \ref{eq:postmean}-\ref{eq:postvar} recursively given the new data points. Furthermore, the approach used for MC sampling of the GPs to be introduced in Section \ref{sec:Gaussian_processes_MC} requires recursive updates of this form as well.

Let the new dataset be given by $\mathcal{D}^+=(\mathcal{D},(\mathbf{z}^+,y^+))$, where $\mathcal{D}$ is the training data for the initial GP, while $\mathbf{z}^{+}$ is the new input and $y^{+}$ the new corresponding output measurement. Then we will refer to the updated mean function and variance function as:
\begin{subequations}
\begin{align}
& \mu_f^+(\mathbf{z};\mathcal{D}^+)    \coloneqq \mathbf{k}^{+ \sf T}(\mathbf{z}) \bm{\Sigma}_{\mathbf{Y}}^{+ -1} \mathbf{Y}^+ \label{eq:postmeanupdate} \\
& \sigma_f^{2+}(\mathbf{z};\mathcal{D}^+) \coloneqq \zeta^2 - \mathbf{k}^{+ \sf T}(\mathbf{z}) \bm{\Sigma}_{\mathbf{Y}}^{+ -1} \mathbf{k}^+(\mathbf{z}) \label{eq:postvarupdate}
\end{align}
\end{subequations}

The updated terms in Equations \ref{eq:postmeanupdate}-\ref{eq:postvarupdate} can be expressed as:
\begin{subequations}
\begin{align} \label{eq:k_update}
    & \mathbf{k}^{+}(\mathbf{z}) = [\mathbf{k}^{\sf T}(\mathbf{z}),k(\mathbf{z},\mathbf{z}^+)]^{\sf T} \\ \label{eq:Y_update} 
    & \mathbf{Y}^+ = [\mathbf{Y}^{\sf T},y^+]^{\sf T}, \quad \mathbf{Z}^+ = [\mathbf{Z}^{\sf T},\mathbf{z}^+]^{\sf T} \\ 
    & \bm{\Sigma}^{+-1}_{\mathbf{Y}} = \begin{bmatrix}
    \bm{\Sigma}_{\mathbf{Y}} & \mathbf{k}^{\sf T}(\mathbf{z}^+)  \\ 
    \mathbf{k}(\mathbf{z}^+) & k(\mathbf{z}^+,\mathbf{z}^+) (+ \sigma_{\nu}^2)
\end{bmatrix}^{-1} \label{eq:update_cov}
\end{align}
\end{subequations}
where $\mathbf{k}(\mathbf{z})$, $\mathbf{Z}$ and $\mathbf{Y}$ refer to quantities of the initial GP. The noise term $\sigma_{\nu}^2$ in the lower diagonal is shown in brackets, since the new "measurement" $y^+$ may be noiseless as is the case for GP MC samples. In this case the noise term should not be added to the new diagonal element.  

Note the updates $\mathbf{k}^{+}(\mathbf{z})$, $\mathbf{Z}^{+}$, and $\mathbf{Y}^+$ are trivial, however the update of the \textit{inverse} covariance matrix $\bm{\Sigma}^{+-1}_{\mathbf{Y}}$ is more involved. In essence we require the inverse of the previous covariance matrix after adding a horizontal row and a vertical row to the covariance matrix of the initial GP, see Equation \ref{eq:update_cov}. For this process there are efficient formula available, one of which is introduced in \ref{app:recursive_update}. These take advantage of the fact that we already know the inverse covariance matrix $\bm{\Sigma}_{\mathbf{Y}}^{-1}$ of the initial GP. Once the update has been carried out, these terms then define the new initial GP. This update procedure is then repeated for the next measurement.  

\subsection{State space model} \label{sec:Gaussian_processes_SS}
In this section we briefly show how the previously introduced GP methodology can be utilized to identify unknown state space models in the form of Equation \ref{eq:f} based on the measurements (data) according to Equation \ref{eq:y}. GPs are commonly applied to model scalar functions with vector inputs as shown in Section \ref{sec:Gaussian_processes_regression}. To extend this to the multi-input, multi-output case as required it is common to build a separate independent GP for each output dimension, see for example \citet{Deisenroth2011}. Let the function in Equation \ref{eq:f} be given by $\mathbf{f}(\mathbf{x},\mathbf{u}) \coloneqq \mathbf{f}(\mathbf{z}) \coloneqq [f_1(\mathbf{z}),\ldots,f_{n_{\mathbf{x}}}(\mathbf{z})]^{\sf T}$. We aim to build a separate GP for each function $f_i(\mathbf{z}) \, \forall i \in \{1,\ldots,n_{\mathbf{x}}\}$ according to Section \ref{sec:Gaussian_processes_regression}. For this purpose we are given observations $\mathbf{Y}_i = [y_{i}^{(1)},\ldots,y_{i}^{(N)}]^{\sf T} \, \forall i \in \{1,\ldots,n_{\mathbf{x}}\}$ and corresponding inputs $\mathbf{Z}=[\mathbf{z}^{(1)},\ldots,\mathbf{z}^{(N)}]^{\sf T}$, where $y_i$ refers to the $i$-th dimension of measurements obtained according to Equation \ref{eq:y}. Let $\mathbf{Y}=[\mathbf{Y}_1,\ldots,\mathbf{Y}_{n_{\mathbf{x}}}]$ correspond to the overall measurements available. The posterior Gaussian distribution of $\mathbf{f}(\cdot)$ at an arbitrary input $\mathbf{z}=(\mathbf{x},\mathbf{u})$ is:  
\begin{subequations} \label{eq:postSS}
\begin{align} 
    & \mathbf{f}(\mathbf{z})|\mathcal{D} \sim \mathcal{N}(\bm{\upmu}_f(\mathbf{z};\mathcal{D}),\bm{\Sigma}_f(\mathbf{z};\mathcal{D})) 
\end{align}
with 
\begin{align} 
    & \bm{\upmu}_f(\mathbf{z};\mathcal{D}) =[\mu_f(\mathbf{z};\mathcal{D}_1),\ldots,\mu_f(\mathbf{z};\mathcal{D}_{n_{\mathbf{x}}})]^{\sf T} \\ 
    & \bm{\Sigma}_f(\mathbf{z};\mathcal{D}) = \diag(\sigma_f^2(\mathbf{z};\mathcal{D}_1),\ldots,\sigma_f^2(\mathbf{z};\mathcal{D}_{n_{\mathbf{x}}})) + \bm{\Sigma}_{\bm{\upomega}}
\end{align}
\end{subequations}
where $\mu_f(\mathbf{z};\mathcal{D}_i)$ and $\sigma_f^2(\mathbf{z};\mathcal{D}_i)$ are the mean function and variance function built according to Section \ref{sec:Gaussian_processes_regression} with datasets $\mathcal{D}_i=(\mathbf{Z},\mathbf{Y}_i) \, \forall i \in \{1,\ldots,n_{\mathbf{x}}\}$ with $\mathcal{D}=(\mathbf{Z},\mathbf{Y})$. 

\begin{remark}[Additive disturbance noise]
Note the additive disturbance noise defined in Equation \ref{eq:f} is simply added to the posterior covariance matrix due to the additive property of multivariate Gaussian distributions.
\end{remark}

In addition, given an initial GP state space model built with a dataset $\mathcal{D}$ and a new data point $(\mathbf{z}^+,\mathbf{y}^+)$, we can update it recursively utilizing the method introduced in Section \ref{sec:Gaussian_processes_online_learning}:
\begin{subequations} \label{eq:updateSS}
\begin{align} \label{eq:mean_ss_update}
    & \bm{\upmu}_f^+(\mathbf{z};\mathcal{D}^+) =[\mu_f^+(\mathbf{z};\mathcal{D}_1^+),\ldots,\mu_f^+(\mathbf{z};\mathcal{D}_{n_{\mathbf{x}}}^+)]^{\sf T} \\
    & \bm{\Sigma}_f^+(\mathbf{z};\mathcal{D}^+) = \diag(\sigma_f^{2+}(\mathbf{z};\mathcal{D}_{1}^+),\ldots,\sigma_f^{2+}(\mathbf{z};\mathcal{D}_{n_{\mathbf{x}}}^+)) + \bm{\Sigma}_{\bm{\upomega}} \label{eq:var_ss_update}
\end{align}
\end{subequations} 
where $\mathcal{D}_i^+=(\mathcal{D}_i,(\mathbf{z}^+,y_i^+)) \, \forall i \in \{1,\ldots,n_{\mathbf{x}}\}$ and $\mathcal{D}^+=(\mathcal{D},(\mathbf{z}^+,\mathbf{y}^+))$

\subsection{Monte Carlo sampling} \label{sec:Gaussian_processes_MC}
GPs are distribution over functions and hence their realizations yield \textit{deterministic} functions, see for example the GP samples shown in Figure \ref{fig:Prior_posterior}. In this section we show how to attain independent samples of GP state space models over a \textit{finite} time horizon. Generating a MC sample of a GP would require sampling an infinite dimensional stochastic process, while there is no known method to achieve this. Instead, approximate approaches have been applied such as spectral sampling \citep{Bradford2018d,QuiCandela2010}. Exact samples of GPs are however possible if the GP MC sample needs to be known at only a finite number of points, which is exactly the situation for state space models over a \textit{finite} time horizon. This technique was first outlined in \citet{Conti2009} and has been employed in \citet{Umlauft2018} for the optimal design of linear controllers. We next outline how to obtain an exact sample of a GP state space model over a \textit{finite} time horizon for an arbitrary feedback control policy.    

Assume we are given a GP state space model as shown in Section \ref{sec:Gaussian_processes_SS} from the input-output dataset $\mathcal{D}=(\mathbf{Z},\mathbf{Y})$. The initial condition $\mathbf{x}_0$ is assumed to follow a known Gaussian distribution as defined in Equation \ref{eq:f}. A GP state space model represents a distribution over possible plant models, for which each realization will lead to a different state sequence. The aim of this section is therefore to show how to obtain a single independent sample of such a state sequence, which can then be repeated to obtain multiple independent MC samples of the GP. Let $\bm{\mathcal{X}}^{(s)}=[\bm{\upchi}_{0}^{(s)},\bm{\upchi}_{1}^{(s)},\ldots,\bm{\upchi}_{T}^{(s)}]^{\sf T}$ denote such a state sequence, where $s$ denotes a particular GP realization and $\bm{\upchi}_{i}^{(s)}$ the realization of the state at discrete time $i$ in the sequence of MC sample $s$. The control inputs at discrete time $i$ for MC sample $s$ are denoted by $\mathpzc{u}_i^{(s)}$. The corresponding joint input of $\bm{\upchi}_{i}^{(s)}$ is represented by $\mathpzc{Z}_i^{(s)} = (\bm{\upchi}_i^{(s)},\mathpzc{u}_i^{(s)})$. We assume the control inputs to be the result of a feedback control policy, which we represent as $\kappa:\mathbb{R}^{n_{\mathbf{x}}} \times \mathbb{R} \rightarrow \mathbb{R}^{n_{\mathbf{u}}}$. The control actions $\mathpzc{u}_i^{(s)}$ are then given as follows:
\begin{equation} \label{eq:control_policy_MC}
    \mathpzc{u}_i^{(s)} = \kappa\left(\bm{\upchi}_{i}^{(s)},i\right)
\end{equation}
where $i$ is the current discrete time. Note the control policy depends on the discrete time directly, since it is a finite horizon control policy.   

Consequently, the control actions over the finite time horizon $T$ are a function of $\mathcal{X}^{(s)}$ and denoted jointly as $\mathcal{U}^{(s)}=[\mathpzc{u}_0^{(s)},\ldots,\bm{\mathpzc{u}}_{T-1}^{(s)}]^{\sf T}=[\kappa\left(\bm{\upchi}_{0}^{(s)},0\right),\ldots,\kappa\left(\bm{\upchi}_{T-1}^{(s)},T-1\right)]^{\sf T}$. Note these control inputs are different for each MC sample $s$ due to feedback.

 We start by sampling the Gaussian distribution of the initial state $\mathbf{x}_0 \sim \mathcal{N}(\bm{\upmu}_{\mathbf{x}_0},\bm{\Sigma}_{\mathbf{x}_0})$ to obtain the realization $\bm{\upchi}_{0}^{(s)}$. The posterior Gaussian distribution of the next state in the sequence $\mathbf{x}_1$ is subsequently given by the GP of $\mathbf{f}(\cdot)$ as defined in Equation \ref{eq:postSS} dependent on $\bm{\upchi}_{0}^{(s)}$:        
\begin{align}
    \mathbf{x}_1 = \mathbf{f}(\mathpzc{Z}_0) \sim \mathcal{N}(\bm{\upmu}_f(\mathpzc{Z}_0;\mathcal{D}),\bm{\Sigma}_f(\mathpzc{Z}_0;\mathcal{D}))
\end{align}

The realization of $\mathbf{x}_1$ is obtained by sampling the above normal distribution, which we will denote as $\bm{\upchi}_{1}^{(s)}$. To obtain the next state in the sequence $\bm{\upchi}_{2}^{(s)}$ we need to first condition on $\bm{\upchi}_{1}^{(s)}$, since this is part of this sampled function path. This requires to treat $\bm{\upchi}_{1}^{(s)}$ similarly to a new training point, however without observation noise (i.e. no $\sigma^2_{\nu}$ is added to the kernel evaluation $k(\mathpzc{Z}_0,\mathpzc{Z}_0)$) and without changing the hyperparameters. Note that if the sampled function were to return to the same input it would lead to the exact same output, since it is conditioned on a noiseless output. This shows that the sampled function is deterministic, since it is the result of sampling.     
Next we draw $\bm{\upchi}_{2}^{(s)}$ according to the posterior Gaussian distribution obtained from adding the previously sampled data point to the training dataset $\mathcal{D}$ as a noiseless observation. This sample is then again added to the training dataset as a noiseless observation, from which the GP is updated and the next state is drawn. This process is repeated until the required time horizon $T$ has been reached. In this paper we consider a finite time horizon control problem and hence the GP state space model should for, moderate time horizons, not become too computationally expensive, since at most $T$ new data-points are added. Nonetheless, for large time horizons this could become a problem and approximate sampling approaches, such as spectral sampling should then be considered instead \citep{Bradford2018d,QuiCandela2010}.   

This sampling approach is summarized in Algorithm 1 below and is illustrated in Figure 2. Each GP MC sample is defined by a state and corresponding control action sequence. Note that this gives us a single MC sample and subsequently needs to be repeated multiple times to obtain multiple realizations. 

\begin{algorithm2e}[H] \label{alg:GPsampling}
 \caption{Gaussian process trajectory sampling}
\Input{$\bm{\upmu}_{\mathbf{x}_0}$, $\bm{\Sigma}_{\mathbf{x}_0}$, $\bm{\upmu}_f(\mathbf{z};\mathcal{D})$, $\bm{\Sigma}_f(\mathbf{z};\mathcal{D})$, $\mathcal{D}$, $T$, $\kappa(\cdot)$}
\Initialize{Draw $\bm{\upchi}_0^{(s)}$ from $\mathbf{x}_0 \sim \mathcal{N}(\bm{\upmu}_{\mathbf{x}_0},\bm{\Sigma}_{\mathbf{x}_0})$.}
 \For{each sampling time $t=1,2,\ldots,T$}{
 \begin{enumerate}
 \item{Determine $\mathpzc{u}_{t-1}^{(s)} = \kappa(\bm{\upchi}_{t-1}^{(s)},t-1)$}
 \item{Draw $\bm{\upchi}_t^{(s)}$ from  $\mathbf{x}_t \sim \mathcal{N}(\bm{\upmu}_f(\mathpzc{Z}_{t-1}^{(s)};\mathcal{D}),\bm{\Sigma}_f(\mathpzc{Z}_{t-1}^{(s)};\mathcal{D}))$, where $\mathpzc{Z}_{t-1}^{(s)} = (\bm{\upchi}_{t-1}^{(s)},\mathpzc{u}_{t-1}^{(s)})$.}
 \item{Define $\mathcal{D}^+ :=(\mathcal{D},(\mathpzc{Z}_{t-1}^{(s)},\bm{\upchi}_t^{(s)}))$, where $\bm{\upchi}_t^{(s)}$ can be viewed as noiseless "measurements".}
  \item{Update the dataset $\mathcal{D}:=([\mathbf{Z}^{\sf T},\mathpzc{Z}_{t-1}^{(s) \sf T}]^{\sf T},[\mathbf{Y},\bm{\upchi}_t^{(s) \sf T}]^{\sf T})$.}
 \item{Recursively update the GP mean function  $\bm{\upmu}_f(\mathbf{z};\mathcal{D}) := \bm{\upmu}^+_f(\mathbf{z};\mathcal{D}^+)$ using Equation \ref{eq:mean_ss_update}.}
 \item{Recursively update the GP covariance function  $\bm{\Sigma}_f(\mathbf{z};\mathcal{D}) := \bm{\Sigma}^+_f(\mathbf{z};\mathcal{D}^+)$ using Equation \ref{eq:var_ss_update}.}
\end{enumerate}}
\Output{State sequence $\bm{\mathcal{X}}^{(s)}= [\bm{\upchi}_0^{(s)},\bm{\upchi}_1^{(s)},\ldots,\bm{\upchi}_T^{(s)}]^{\sf T}$ and control sequence $\bm{\mathcal{U}}^{(s)}=[\mathpzc{u}_0^{(s)},\ldots,\bm{\mathpzc{u}}_{T-1}^{(s)}]^{\sf T}$}
\end{algorithm2e}

\begin{figure}[H] \centering
   \includegraphics[width=0.8\textwidth]{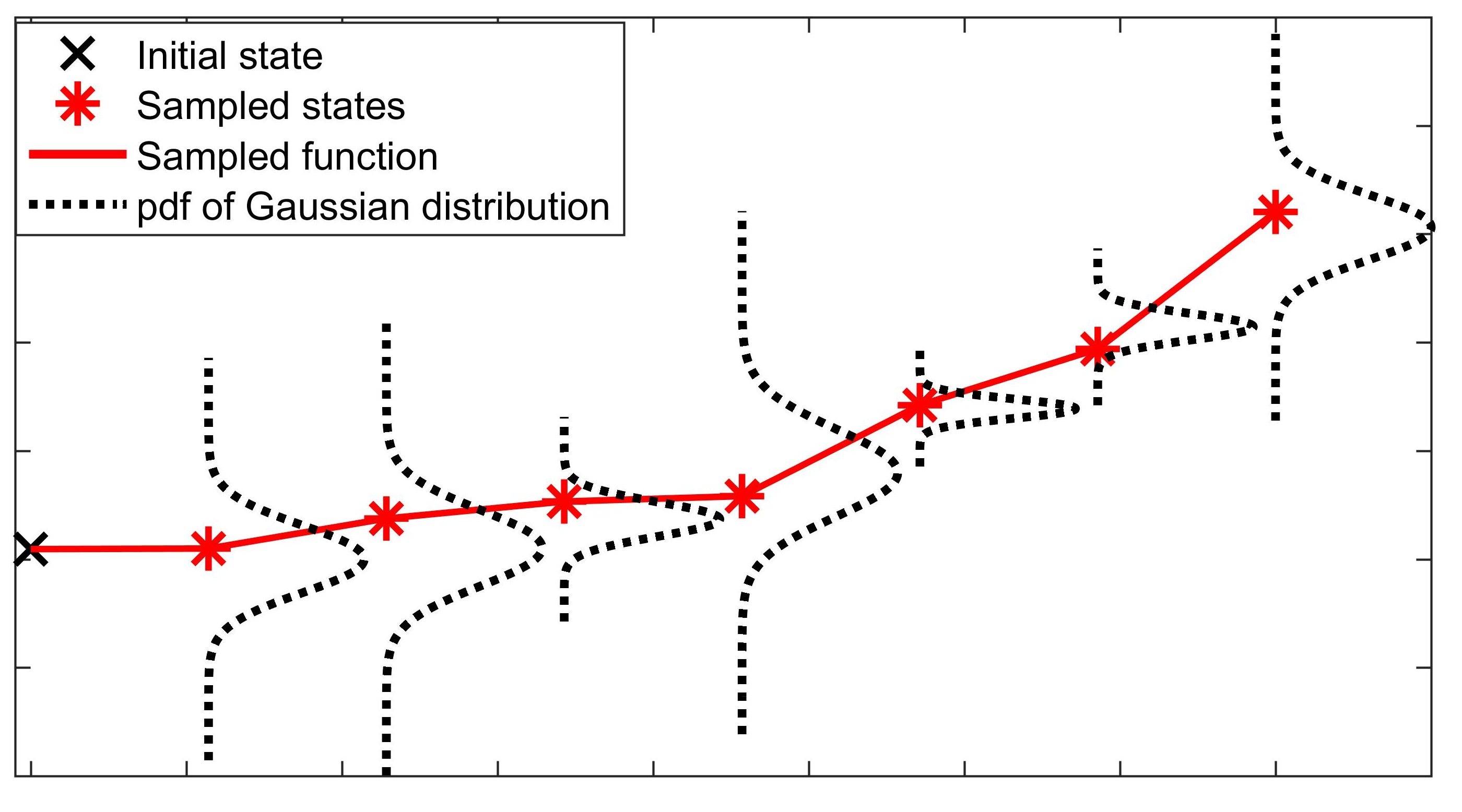}
  \caption{Illustration of GP sampling scheme for a 1 dimensional function.}
  \label{fig:GP_sampling}
\end{figure}

Lastly, we define a \textit{nominal} trajectory by setting all samples in Algorithm 1 to their mean values. Let $\overline{\bm{\mathcal{X}}}= [\overline{\bm{\upchi}}_0,\overline{\bm{\upchi}}_1,\ldots,\overline{\bm{\upchi}}_T]^{\sf T}$ refer to this \textit{nominal} state sequence and  $\overline{\bm{\mathcal{U}}}=[\overline{\mathpzc{u}}_0,\ldots,\overline{\mathpzc{u}}_{T-1}]^{\sf T}$ to the corresponding \textit{nominal} control sequence, where $\overline{\bm{\upchi}}_i$ and $\overline{\mathpzc{u}}_i$ are the values of the states and control inputs of the \textit{nominal} trajectory at discrete time $i$ respectively. Therefore by definition, $\overline{\bm{\upchi}}_0=\bm{\upmu}_{\mathbf{x}_0}$, $\overline{\mathpzc{u}}_{t-1}=\kappa(\overline{\bm{\upchi}}_{t-1},t-1)$, and $\overline{\bm{\upchi}}_t=\bm{\upmu}_f(\overline{\mathpzc{Z}}_{t-1};\mathcal{D})$, where $\overline{\mathpzc{Z}}_{t-1}=(\overline{\bm{\upchi}}_{t-1},\overline{\mathpzc{u}}_{t-1})$. Note updating $\bm{\upmu}_f(\cdot;\mathcal{D})$ with mean values has no effect.    

In Section \ref{sec:Solution_approach} a NMPC formulation is introduced, which uses the mean function of the GP as the prediction model. The control actions from this NMPC algorithm then define the control policy in Equation \ref{eq:control_policy_MC}, while the MC samples represent possible plant responses. This is then exploited to tune the GP NMPC algorithm to attain the desired behaviour.   

\section{Solution approach} \label{sec:Solution_approach}
From the input-output dataset $\mathcal{D}=(\mathbf{Z},\mathbf{Y})$ defined in Equation \ref{eq:datasetdef}, we fit a GP state space model as outlined in Section \ref{sec:Gaussian_processes_SS}. The aim is to solve the problem defined in Section \ref{sec:prob_def} based on this GP state space model. In this context the GP represents a distribution over possible plant models for the process given the available dataset. In Section \ref{sec:Gaussian_processes_MC} we have shown how to create a sample of this plant model over a finite time horizon $T$, which each lead to different state sequences and corresponding control sequences based on a control policy. In this paper we aim to design a NMPC algorithm based on the GP that acts as this control policy. The MC samples are utilized to tune the NMPC formulation by adjusting so-called back-offs to tighten the constraints and attain the closed-loop probabilistic constraint satisfaction. We next state the GP NMPC formulation, which is based on the tightened constraint set and predictions from the mean function $\bm{\upmu}_f(\cdot;\mathcal{D})$ and covariance function $\bm{\Sigma}_f(\cdot;\mathcal{D})$.

\subsection{Finite-horizon Gaussian process model predictive control formulation}
In this section we define the NMPC OCP based on the GP \textit{nominal} model given the dataset $\mathcal{D}=(\mathbf{Z},\mathbf{Y})$. For the GP NMPC formulation the initial state $\mathbf{x}$ at each sampling time is assumed to be measured or estimated and propagated forward in time exploiting the GP mean function. The predicted states are then used to optimize the objective subject to the tightened constraints. Let the corresponding optimization problem be denoted as $P_T\left(\bm{\upmu}_f(\cdot;\mathcal{D}),\bm{\Sigma}_f(\cdot;\mathcal{D});\mathbf{x},t\right)$ for the current \textit{known} state $\mathbf{x}$ at discrete time $t$ based on the mean function $\bm{\upmu}_f(\cdot;\mathcal{D})$ and covariance function $\bm{\Sigma}_f(\cdot;\mathcal{D})$:
\begin{equation}
\begin{aligned} \label{eq:nominalMPC}
& \underset{\hat{\mathbf{U}}_{t:T-1}}{\text{minimize}} \quad  \hat{V}_T(\mathbf{x},t,\hat{\mathbf{U}}_{t:T-1}) = \sum_{k=t+1}^{T-1} \left[ \ell(\hat{\mathbf{x}}_k,\hat{\mathbf{u}}_k) + \eta_k \tr\left(\bm{\Sigma}_f(\hat{\mathbf{x}}_k;\mathcal{D})\right) \right] +  \ell_f(\hat{\mathbf{x}}_T) \\ 
& \text{subject to:}  \\
& \hat{\mathbf{x}}_{k+1} = \bm{\upmu}_f(\hat{\mathbf{z}}_k;\mathcal{D}), \quad \hat{\mathbf{z}}_k = (\hat{\mathbf{x}}_k,\hat{\mathbf{u}}_k)   && \forall k \in \{t,\ldots,T-1\} \\
& \hat{\mathbf{x}}_{k+1} \in  \overline{\mathbb{X}}_{k+1}, \quad \hat{\mathbf{u}}_k \in \mathbb{U}_k && \forall k \in \{t,\ldots,T-1\} \\
& \hat{\mathbf{x}}_t = \mathbf{x}
\end{aligned}
\end{equation}
where $\hat{\mathbf{x}}$, $\hat{\mathbf{u}}$, and $\hat{V}_T(\cdot)$ refers to the states, control inputs, and control objective of the MPC formulation, $\hat{\mathbf{U}}_{t:T-1} = [\hat{\mathbf{u}}_t,\ldots,\hat{\mathbf{u}}_{T-1}]^{\sf T}$, $\eta_k$ are weighting factors to penalize uncertainty, and $\overline{\mathbb{X}}_k$ is a tightened constraint set denoted by: $\overline{\mathbb{X}}_{k} = \{\mathbf{x} \in \mathbb{R}^{n_x} \; | \; g_j^{(k)}(\mathbf{x}) + b_j^{(k)} \leq 0, \, j=1,\ldots,n_g \}$. The variables $b_j^{(k)}$ represent so-called back-offs, which tighten the original constraints $\mathbb{X}_t$ defined in Equation \ref{eq:xcon}.

\begin{remark}[Objective in expectation]
It should be noted that the above objective in Equation \ref{eq:nominalMPC} does not exactly optimize the objective in Equation \ref{eq:objdef}, since it is difficult to obtain the expectation of a nonlinear function. Approximations of this can be found in \citet{Hewing2017}, however these generally are considerably more expensive and often only lead to marginally improved performance.
\end{remark}

\begin{remark}[Scaling for state dependency factors]
Note we have opted for scalar scaling factors $\eta_k$ to account for state dependency. For this to work reliably it is therefore necessary to normalize the data to ensure that all data has approximately the same magnitude, for example normalizing the data to have zero mean and unit variance. 
\end{remark}

The NMPC algorithm solves $P_T\left(\bm{\upmu}_f(\cdot;\mathcal{D}),\bm{\Sigma}_f(\cdot;\mathcal{D});\mathbf{x}_t,t\right)$ at each sampling time $t$ given the current state $\mathbf{x}_t$ to obtain an optimal control sequence:
\begin{align}
\hat{\mathbf{U}}^*_{t:T-1}\left(\bm{\upmu}_f(\cdot;\mathcal{D}),\bm{\Sigma}_f(\cdot;\mathcal{D});\mathbf{x}_t,t\right) = [\hat{\mathbf{u}}^*_t\left(\bm{\upmu}_f(\cdot;\mathcal{D}),\bm{\Sigma}_f(\cdot;\mathcal{D});\mathbf{x}_t,t\right),\ldots,\hat{\mathbf{u}}^*_{T-1}\left(\bm{\upmu}_f(\cdot;\mathcal{D}),\bm{\Sigma}_f(\cdot;\mathcal{D});\mathbf{x}_t,t\right)]^{\sf T} 
\end{align}

Only the first optimal control action is applied to the plant at time $t$ before the same optimization problem is solved at time $t+1$ with a new state measurement $\mathbf{x}_{t+1}$. This procedure implicitly defines the following feedback control law, which needs to be repeatedly solved for each new measurement $\mathbf{x}_t$:
\begin{align} \label{eq:GP_control_policy}
\kappa(\bm{\upmu}_f(\cdot;\mathcal{D}),\bm{\Sigma}_f(\cdot;\mathcal{D});\mathbf{x}_t,t) = \hat{\mathbf{u}}^*_t\left(\bm{\upmu}_f(\cdot;\mathcal{D}),\bm{\Sigma}_f(\cdot;\mathcal{D});\mathbf{x}_t,t\right)
\end{align}

It is explicitly denoted that the control actions depend on the GP model used. There are several important variations of the GP NMPC control policy. Firstly, it may seem reasonable to update the mean and covariance function using the previous state measurement and corresponding input by applying the recursive update rules introduced in Section \ref{sec:Gaussian_processes_online_learning}. We will refer to this as \textit{learning}. This may however lead to a more expensive and less reliable NMPC algorithm. 

Secondly, the algorithm may want to avoid regions in which there is great uncertainty due to sparsity of data. This can be achieved by assigning some of the $\eta_k$ with non-zero values to penalize the algorithm moving into these regions with high variance. This explicitly takes advantage of the state dependency of the noise covariance function and is hence referred as \textit{state dependent}. It should be noted that evaluation of the covariance function is computationally expensive. It was determined that setting only $\eta_{t+1}$ to a non-zero value is often sufficient due to the continued feedback update, i.e. penalizing only the variance for the one-step ahead prediction at time $t$. These variations can be summarized as follows:
\begin{itemize}
    \item{\textit{Learning}: Update the mean and variance function of the GP using the previous state measurement and the known corresponding input.}
    \item{\textit{State dependent}: Set some $\eta_k$ not equal to zero, which will lead to the NMPC algorithm trying to find a path that has less variance and hence exploiting the state dependent nature of the uncertainty.}
    \item{\textit{Non learning}: Keep the mean and variance function the same throughout the run.}
    \item{\textit{Non state dependent}: Set all $\eta_k$ to zero and hence ignoring the possible state dependency of the uncertainty.}
\end{itemize}

\begin{remark}[Full state feedback]
Note in the control algorithm we have assumed full state feedback, i.e. it is assumed that the full state can be measured without noise. This assumption can be dropped if required by introducing a suitable observer and introduced in the closed-loop simulations to account for this additional uncertainty.  
\end{remark}

\subsection{Probabilistic guarantees}
In this section we illustrate how to obtain probabilistic guarantees for the joint chance constraint introduced in Section \ref{sec:prob_def} in Equation \ref{eq:xcon} based on independent samples of the GP plant model. For convenience we define a single-variate random variable $C(\cdot)$ representing the satisfaction of the joint chance constraint \citep{Curtis2018}:
\begin{subequations} \label{eq:c_def}
\begin{align}
    & C(\mathbf{X}) = \inf_{(j,t) \in \{1,\ldots,n_g\} \times \{0,\ldots,T\}} {g_j^{(t)}(\mathbf{x}_t)} \\
    & \mathbb{P}\left\{C(\mathbf{X}) \leq 0\right\} = \mathbb{P} \left\{ \bigcap^T_{t=0} \{ \mathbf{x}_t \in \mathbb{X}_t \}    \right\}
\end{align}
\end{subequations}
where $\mathbf{X} = [\mathbf{x}_0,\ldots,\mathbf{x}_T]^{\sf T}$ defines a state sequence, and $\mathbb{X}_t = \{ \mathbf{x} \in \mathbb{R}^{n_{\mathbf{x}}} \mid g_j^{(t)}(\mathbf{x}) \leq 0, j=1,\ldots,n_g \}$. 

The probability in Equation \ref{eq:c_def} is intractable, however a good nonparametric approximation is often achieved utilizing the so-called empirical cumulative distribution function (ecdf). We define the cdf to be approximated as follows:
\begin{align} \label{eq:cdf_c}
    F_{C(\mathbf{X})}(c) = \mathbb{P}\{C(\mathbf{X}) \leq c \} 
\end{align}

Assuming we are given $S$ independent and identically distributed MC samples of $\mathbf{X}$ and hence of $C(\mathbf{X})$, the ecdf estimate of the true cdf in Equation \ref{eq:cdf_c} is given by:
\begin{equation} \label{eq:approx_joint}
    F_{C(\mathbf{X})}(c) \approx \hat{F}_{C(\mathbf{X})}(c) = \frac{1}{S} \sum_{s=1}^{S} \mathbf{1}\{C(\bm{\mathcal{X}}^{(s)}) \leq c\} 
\end{equation}
where $\bm{\mathcal{X}}^{(s)}$ is the $s$-th MC sample and $\hat{F}_{C(\mathbf{X})}(c)$ is the ecdf approximation of the true cdf $F_{C(\mathbf{X})}(c)$. 

The quality of the approximation in Equation \ref{eq:approx_joint} strongly depends on the number of samples used and it is therefore desirable to quantify the residual uncertainty of the sample approximation. This problem has been studied to a great extent in the statistics literature \citep{Clopper1934}. In addition, there are several works applying these results for chance constrained optimization, see for example \citet{Alamo2015}. The main result applied in this study is given below in Theorem 1.  

\begin{figure}[H] \centering
   \includegraphics[width=0.75\textwidth]{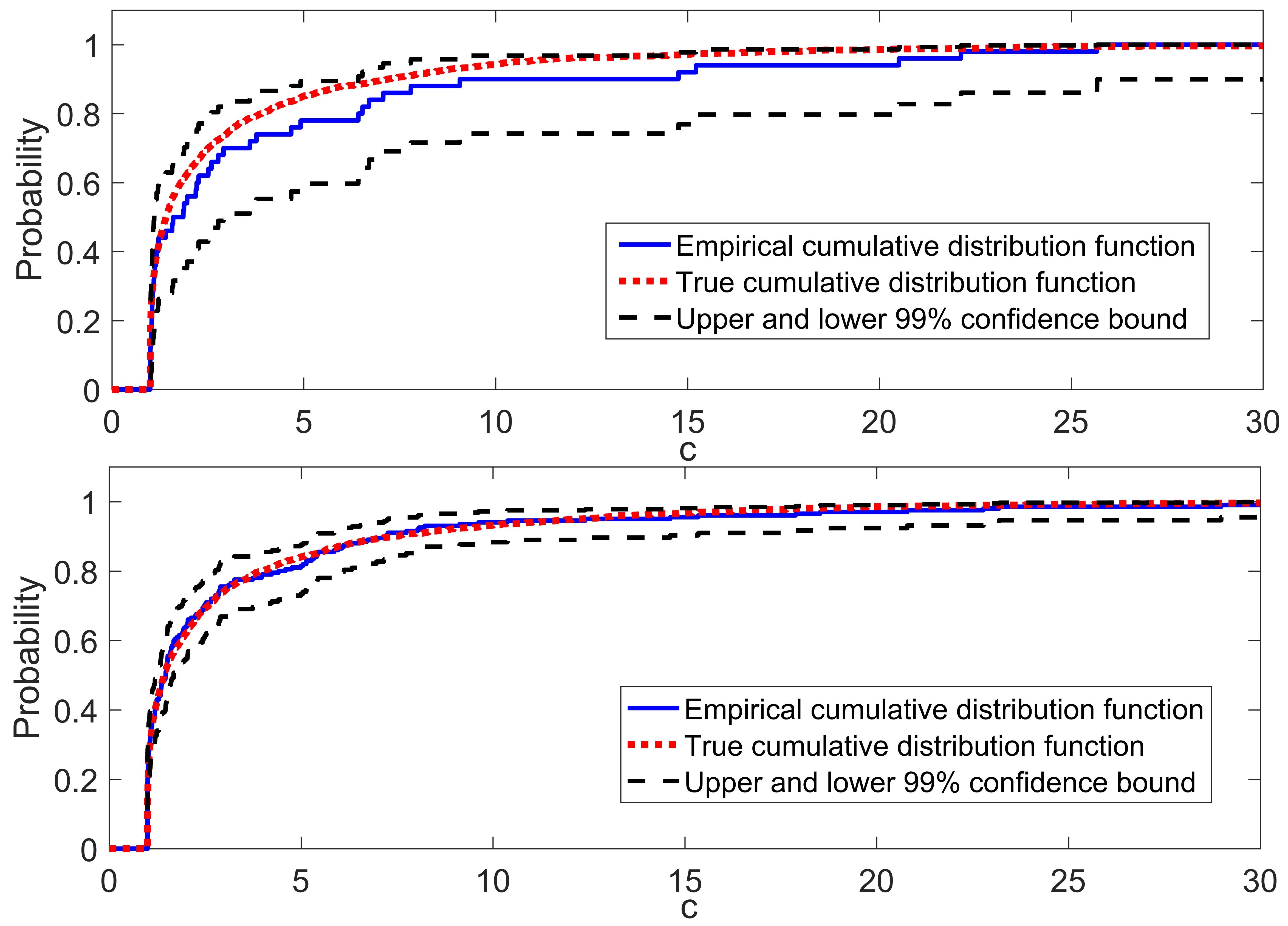}
  \caption{Illustration of the cdf confidence bound of $\alpha=0.01$ for sample sizes of $S=50$ (top) and $S=200$ (bottom).}
  \label{fig:CDF_confidence}
\end{figure}

\begin{theorem}[Confidence interval for empirical cumulative distribution function] Assume we are given a value of the ecdf, $\hat{\beta} = \hat{F}_{C(\mathbf{X})}(c)$, as defined in Equation \ref{eq:approx_joint} based on $S$ independent samples of $C(\mathbf{X})$, then the true value of the cdf, $\beta = F_{C(\mathbf{X})}(c)$, as defined in Equation \ref{eq:cdf_c} has the following lower and upper confidence bounds:
\begin{subequations}
\begin{align}
    & \mathbb{P}\left\{ \beta \geq \hat{\beta}_{lb} \right\} \geq 1-\alpha, && \hat{\beta}_{lb} = \textup{betainv}\left(\alpha, S \hat{\beta}, S - S \hat{\beta} + 1 \right),   \\
    & \mathbb{P}\left\{ \beta \leq \hat{\beta}_{ub} \right\} \geq 1-\alpha, && \hat{\beta}_{ub} = \textup{betainv}\left(1-\alpha, S \hat{\beta} + 1, S - S \hat{\beta} \right).
\end{align}
\end{subequations}
\end{theorem}
\begin{proof}
The proof uses standard results in statistics and can be found in \citet{Clopper1934,Streif2014}. The proof relies on the following observations. Firstly, $\mathbf{1}\{C(\bm{\mathcal{X}}^{(s)}) \leq c\}$ for a \textit{fixed value} of $c$ describes a Bernoulli random variable, in which either $C(\bm{\mathcal{X}}^{(s)})$ exceeds $c$ and takes the value $0$ or otherwise takes the value $1$ with probability $F_{C(\mathbf{X})}(c)$. Secondly, the ecdf describes the number of \textit{successes} of $S$ realizations of these Bernoulli random variables divided by the total number of samples and hence follows a binomial distribution, i.e. $\hat{F}_{C(\mathbf{X})}(c) \sim \frac{1}{N_S}\text{Bin}(S,F_{C(\mathbf{X})}(c))$. The confidence bound for the ecdf can consequently be determined from the Binomial cdf. This method was first introduced by \citet{Clopper1934} as "\textit{exact confidence intervals}". Due to the close relationship between beta distributions and binomial distributions, a simplified expression can be obtained using beta distributions instead, which leads to the theorem shown \citep{Streif2014}. 
\end{proof}
In other words the probability of $\beta$ exceeding the value $\hat{\beta}_{ub}$ has a probability of $\alpha$ and the probability of $\beta$ being less than or equal to $\hat{\beta}_{lb}$ has also a probability of $\alpha$. In particular, $\hat{\beta}_{lb}$ for small $\alpha$ represents a conservative lower bound on the true probability $\beta$. An illustration of the confidence bound for the ecdf is shown in Figure \ref{fig:CDF_confidence}. In general more samples will lead to a tighter confidence bound as expected. The theorem provides a lower bound $\hat{\beta}_{lb}$ that accounts for the statistical error due to the finite sample estimate made, i.e. it gives us a conservative value that is less than or equal to the true probability of feasibility with a confidence level of $1-\alpha$.  

We assume we are given $S$ \textit{independent} samples of the trajectory $\mathbf{X}$ generated according to Section \ref{sec:Gaussian_processes_MC}. Let the approximate ecdf be given by $\hat{\beta}=\hat{F}_{C(\mathbf{X})}(0)$ according to Equation \ref{eq:approx_joint}, and $\hat{\beta}_{lb}$ the corresponding lower bound according to Theorem 1 with confidence level $1-\alpha$. From this the following Corollary follows:

\begin{corollary}[Feasibility probability]
Assuming the GP representation of the plant model to be a correct description of the uncertainty of the system and given a value of the ecdf $\hat{\beta} = \hat{F}_{C(\mathbf{X})}(0)$, as defined in Equation \ref{eq:approx_joint} based on $S$ independent samples and a corresponding lower bound $\hat{\beta}_{lb} \geq 1-\epsilon$ with a confidence level of $1-\alpha$, then the original chance constraint in Equation \ref{eq:xcon} holds true with a probability of at least $1-\alpha$.  
\end{corollary}
\begin{proof}
The GP MC sample described in Section \ref{sec:Gaussian_processes_MC} is exact and therefore each sample of the GP plant state space model leads to independent state trajectories $\mathbf{X}$ according to the GP distribution. From $S$ such samples a valid lower bound $\hat{\beta}_{lb}$ to the true cdf value $\beta$ can be determined from Theorem 1 with a confidence level of $1-\alpha$. If $\hat{\beta}_{lb}$ is greater than or equal to $1-\epsilon$, then the following probabilistic bound holds on the true cdf value $\beta$ according to Theorem 1: $\mathbb{P}\left\{\beta \geq \hat{\beta}_{lb} \geq  1 - \epsilon \right\} \geq 1-\alpha$, which in other words means that $\beta = \mathbb{P}\left\{ C(\mathbf{X}) \leq 0 \right\} \geq 1-\epsilon$ with a probability of at least $1-\alpha$. 
\end{proof}

\subsection{Determining back-off constraints} \label{sec:back_off_iterations}
In this section we describe how to determine the required back-off values to tighten the constraints of the GP NMPC algorithm in Equation \ref{eq:nominalMPC}. The aim is to choose these values to obtain probabilistic guarantees for the state constraints defined in Equation \ref{eq:xcon} despite not knowing the exact dynamics. The GP provides a nominal model for the GP NMPC formulation in Equation \ref{eq:nominalMPC} using the mean function and a distribution of possible plant models given the initial dataset in Equation \ref{eq:datasetdef}. This distribution is exploited to attain different possible plant model realizations to simulate the closed-loop response of the GP NMPC algorithm and then uses the response values to tighten the constraints, such that the original constraint set is satisfied with a high probability according to Equation \ref{eq:xcon}. 

We have shown how to obtain the closed-loop trajectory of a control policy according to the realization of a plant model GP distribution in Section \ref{sec:Gaussian_processes_MC}. To this we apply the GP NMPC control policy defined in Equation \ref{eq:GP_control_policy}. We propose to utilize $S$ independent MC samples of the GP distribution generated according to Section \ref{sec:Gaussian_processes_MC}, which then in turn describe $S$ different possible plant models with corresponding state and control trajectories. The goal now is to adjust the back-offs, such that the $S$ different state trajectories adhere the original constraint set for all but a few samples to attain the required probability of constraint satisfaction. Let $\bm{\mathcal{X}}^{(s)}=[\bm{\upchi}_0^{(s)},\ldots,\bm{\upchi}_T^{(s)}]^{\sf T}$ refer to the state trajectory of sample $s$. 

The update rule for adjusting the back-offs is based on two steps: First, we define an approximate constraint set, which is then adjusted by a constant factor to obtain the required constraint satisfaction using the ecdf of the joint chance constraint in Equation \ref{eq:ecdfg}. The approximate constraint set in essence needs to reflect the difference in constraint values for the nominal model of the MPC and the realizations of the GP plant model. We first set all the back-off values to zero and run $S$ MC samples of the GP plant model. As defined in Section \ref{sec:Gaussian_processes_MC}, let $\overline{\bm{\upchi}}_t$ refer to the states according to the \textit{nominal} trajectory with the back-offs set to zero as well. Now assume we aim to obtain back-off values that imply the satisfaction of the following individual chance constraints:
\begin{align} \label{eq:chance_implication}
    & g_j^{(t)}(\overline{\bm{\upchi}}_t) + b_j^{(t)} = 0 \implies \mathbb{P}\left\{g_j^{(t)}(\bm{\upchi}_t) \leq 0 \right\} \geq 1 - \delta 
\end{align}
where $\delta$ is a tuning parameter and should be set to a reasonably low value. 

The rule in other words aims to find approximate back-offs for the nominal predictions $\overline{\bm{\upchi}}_t$ utilized in the MPC formulation in Equation \ref{eq:nominalMPC}, such that the chance constraints holds for any possible GP plant model MC sample with a probability of $1-\delta$. The parameter $\delta$ in this case is a tuning parameter to obtain the initial back-off values. Note the considered individual chance constraint in Equation \ref{eq:chance_implication} is only there to obtain the \textit{initial} constraint set and is unrelated to the joint chance constraint in Equation \ref{eq:xcon}, which we fulfil by adjusting this approximate constraint set using a root-finding algorithm.    

To accomplish this we define the following ecdf based on the $S$ MC samples available:
\begin{align} \label{eq:ecdfg}
    & \hat{F}_{g_j^{(t)}}(0) = \frac{1}{S} \sum_{s=1}^S \mathbf{1}\{g_j^{(t)}({\bm{\upchi}_t^{(s)}}) \leq 0 \}
\end{align}
where $\hat{F}_{g_j^{(t)}}(0)$ is a sample approximation of the chance constraint given in Equation (\ref{eq:chance_implication}) on the RHS.

In \citet{Paulson2018} it is proposed to employ the inverse ecdf to approximately fulfill the requirement given in Equation (\ref{eq:chance_implication}) using the $S$ MC samples available. The back-offs can then be stated as follows:
\begin{equation} \label{eq:update_rule}
\tilde{b}_j^{(t)} = \hat{F}_{g_j^{(t)}}^{-1}(1-\delta) - g_j^{(t)}(\overline{\bm{\upchi}}_t) \quad \forall (j,t) \in \{1,\ldots,n_g^{(t)}\} \times \{1,\ldots,T\} 
\end{equation}
where $\hat{F}_{g_j^{(t)}}^{-1}$ is the inverse of the ecdf defined in Equation \ref{eq:ecdfg} and $\tilde{b}_j^{(t)}$ refers to these initial back-off values. Note the inverse of an ecdf is given by the quantile function of the discrete $S$ values with probability $\delta$.     
This gives us the initial back-off values as required for the first step. Note that both the nominal trajectory $\overline{\bm{\upchi}}_t$ and the GP realizations depend on the back-off values, which is however ignored since we are only interested in obtaining some reasonable initial values. 
In the next step these back-off values are further adjusted using a constant $\textit{back-off}$ factor $\gamma$. The new back-offs are then defined as:
\begin{equation}
   b_j^{(t)} = \gamma \tilde{b}_j^{(t)} \quad \forall (j,t) \in \{1,\ldots,n_g^{(t)}\} \times \{1,\ldots,T\}
\end{equation}

We aim to change $\gamma$ until the lower bound of the ecdf $\hat{\beta}_{lb}$ as defined in the previous section for the joint chance constraint is equal to $1-\epsilon$ in Equation (\ref{eq:xcon}). This is a root finding problem, in which $\gamma$ is adjusted until $\hat{\beta}_{lb}$ reaches the required value:
\begin{equation}
    h(\gamma) = \hat{\beta}_{lb}(\gamma) - (1 - \epsilon)
\end{equation}
where the aim is to determine a value of $\gamma$, such that $h(\gamma)$ is approximately zero. 

To attain the required $\gamma$ we use the so-called \textit{bisection method} \citep{Beers2007}. This method determines the root of a function in an interval $a_{\gamma}$ and $b_{\gamma}$, where $h(a_{\gamma})$ and $h(b_{\gamma})$ have opposite signs.  In our case this is relatively easy. Setting the value of $\gamma$ too low returns generally a negative value of $h(\gamma)$ due to the constraint violations using low back-offs, while setting it too high leads to positive values leading to a too conservative solution. Note we generally set the initial $a_{\gamma}$ to zero since this corresponds to the $S$ MC samples used to determine $\tilde{b}_j^{(t)}$. The \textit{bisection method} consists of repeatably bisecting the interval, in which the root is contained. The overall algorithm to determine the required back-offs in $n_b$ back-off iterations is summarized below as Algorithm 2. The output of the algorithm are the required back-offs with the corresponding lower bound on the probability of satisfying the state chance constraint. Note for $\textit{learning}=true$ the mean and covariance function of the GP are recursively updated utilizing the same procedure as for the update of the GP plant model MC sample.

\begin{remark}[Conservativeness of chance constraint]
Note to adjust the back-offs we use the ecdf, which does account for the \textit{true} shape of the underlying probability distribution. This avoids the problem that is often faced in stochastic optimization utilizing Chebyshev's inequality to robustly approximate chance constraints, which is often excessively conservative \citep{Paulson}.  
\end{remark}

\begin{algorithm2e}[H] \label{alg:back_off_algorithm}
 \caption{Back-off iterative updates}
\Input{$\bm{\upmu}_{\mathbf{x}_0}$, $\bm{\Sigma}_{\mathbf{x}_0}$, $\bm{\upmu}_f(\mathbf{z};\mathcal{D})$, $\bm{\Sigma}_f(\mathbf{z};\mathcal{D})$, $\mathcal{D}$, $T$, $V_T(\mathbf{x},t,\hat{\mathbf{U}}_{t:T-1})$, $\mathbb{X}_t$, $\mathbb{U}_t$, $\epsilon$, $\alpha$, $\delta$, learning, S, $n_b$}

\Initialize{Set all $b_j^{(t)}=0$ and $\delta$ to some reasonable value, set $a_{\gamma}=0$ and $b_{\gamma}$ to some reasonably high value, such that $b_{\gamma} - (1-\epsilon)$ has a positive sign. Define $\mathcal{D}_0 = \mathcal{D}$ as the initial dataset.}

\For{$n_b$ back-off iterations}{
\uIf{$n_b > 0$}{$c_{\gamma} := (a_{\gamma} + b_{\gamma})/2$ \\ $b_j^{(t)} := c_{\gamma} \tilde{b}_j^{(t)} \quad (j,t) \in \{1,\ldots,n_g^{(t)}\} \times \{1,\ldots,T\}$} 
\For{each MC sample $s=1,2,\ldots,S$}{
$\mathcal{D} := \mathcal{D}_0$ \\
Draw $\bm{\upchi}_0^{(s)}$ from $\mathbf{x}_0 \sim \mathcal{N}(\bm{\upmu}_{\mathbf{x}_0},\bm{\Sigma}_{\mathbf{x}_0})$ \\
 \For{each sampling time $t=1,2,\ldots,T$}{
 \uIf{$learning = true$}{$\text{1.}$ Determine $\mathpzc{u}_{t-1}^{(s)} = \kappa(\bm{\upmu}_f(\cdot;\mathcal{D}),\bm{\Sigma}_f(\cdot;\mathcal{D});\mathbf{x}_t,t)$.}
 
 \Else{$\text{1.}$ Determine $\mathpzc{u}_{t-1}^{(s)} = \kappa(\bm{\upmu}_f(\cdot;\mathcal{D}_0),\bm{\Sigma}_f(\cdot;\mathcal{D}_0);\mathbf{x}_t,t)$.}

 \begin{enumerate} \setcounter{enumi}{1}
 \item{Draw $\bm{\upchi}_t^{(s)}$ from  $\mathbf{x}_t \sim \mathcal{N}(\bm{\upmu}_f(\mathpzc{Z}_{t-1}^{(s)};\mathcal{D}),\bm{\Sigma}_f(\mathpzc{Z}_{t-1}^{(s)};\mathcal{D}))$, where $\mathpzc{Z}_{t-1}^{(s)} = (\bm{\upchi}_{t-1}^{(s)},\mathpzc{u}_{t-1}^{(s)})$.}
 \item{Define $\mathcal{D}^+ :=(\mathcal{D},(\mathpzc{Z}_{t-1}^{(s)},\bm{\upchi}_t^{(s)}))$, where $\bm{\upchi}_t^{(s)}$ can be viewed as noiseless "measurements".}
  \item{Update the dataset $\mathcal{D}:=([\mathbf{Z}^{\sf T},\mathpzc{Z}_{t-1}^{(s) \sf T}]^{\sf T},[\mathbf{Y},\bm{\upchi}_t^{(s) \sf T}]^{\sf T})$.}
 \item{Recursively update the GP mean function $\bm{\upmu}_f(\mathbf{z};\mathcal{D}) := \bm{\upmu}^+_f(\mathbf{z};\mathcal{D}^+)$ using Equation \ref{eq:mean_ss_update}.}
 \item{Recursively update the GP covariance function $\bm{\Sigma}_f(\mathbf{z};\mathcal{D}) := \bm{\Sigma}^+_f(\mathbf{z};\mathcal{D}^+)$ using Equation \ref{eq:var_ss_update}.}
 \item{Define $\bm{\mathcal{X}}^{(s)}= [\bm{\upchi}_0^{(s)},\bm{\upchi}_1^{(s)},\ldots,\bm{\upchi}_T^{(s)}]^{\sf T}$ and $\bm{\mathcal{U}}^{(s)}=[\mathpzc{u}_0^{(s)},\ldots,\bm{\mathpzc{u}}_{T-1}^{(s)}]^{\sf T}$.}
\end{enumerate}}}

$\hat{\beta} := \hat{F}_{C(\bm{\mathcal{X}}^{(s)})}(0) = \frac{1}{S} \sum_{s=1}^{S} \mathbf{1}\{C(\bm{\mathcal{X}}^{(s)}) \leq 0\}$  \\
$\hat{\beta}_{lb} := 1-\text{betainv}\left(\alpha,S + 1 - S \hat{\beta}, S \hat{\beta} \right)$ 

\uIf{$nb = 0$}{
Let $\tilde{b}_j^{(t)} = \hat{F}_{g_j^{(t)}}^{-1}(\delta) - g_j^{(t)}(\overline{\bm{\upchi}}_t) \quad \forall (j,t) \in \{1,\ldots,n_g^{(t)}\} \times \{1,\ldots,T\}$ \\
$\hat{\beta}_{lb}^{a_{\gamma}} := \hat{\beta}_{lb} - (1-\epsilon) $} 

\Else{$\hat{\beta}_{lb}^{c_{\gamma}} := \hat{\beta}_{lb} - (1-\epsilon)$ \\ 
\uIf{ $\text{sign}(\hat{\beta}_{lb}^{c_{\gamma}}) = \text{sign}(\hat{\beta}_{lb}^{a_{\gamma}})$}{$a_{\gamma} := c_{\gamma}$ \\
$\hat{\beta}_{lb}^{a_{\gamma}} := \hat{\beta}_{lb}^{c_{\gamma}}$}
\Else{$b_{\gamma} := c_{\gamma}$}}}

\Output{$b_j^{(t)} \quad \forall (j,t) \in \{1,\ldots,n_g^{(t)}\} \times \{1,\ldots,T\}, \hat{\beta}_{lb}$}
\end{algorithm2e}

\subsection{Algorithm}
The overall algorithm proposed in this paper is summarized in this section. Firstly, the problem to be solved needs to be defined as outlined in Section \ref{sec:prob_def}. Thereafter, it needs to be decided if the GP NMPC should \textit{learn} online or exploit the state dependency of the uncertainty as shown in Section \ref{sec:Solution_approach}. Once the GP NMPC has been formulated the back-offs are determined by running closed-loop simulations of the defined problem as shown in Section \ref{sec:back_off_iterations}. Lastly, these back-offs then give us the tightened constraint set required for the GP NMPC \textit{online}. This GP NMPC is then run online solving the problem initially outlined. An overall summary can be found in Algorithm \ref{alg:algorithm_summary}. 

\begin{algorithm2e}[H] \label{alg:algorithm_summary}
 \caption{Back-off GP NMPC}
 \textit{Offline Computations}
 \begin{enumerate}
\item{Build GP state-space model from data-set $\mathcal{D}=(\mathbf{Z},\mathbf{Y})$ and additive disturbance $\bm{\Sigma}_{\bm{\upomega}}$.}
\item{Choose time horizon $T$, initial condition mean $\bm{\upmu}_{\mathbf{x}_0}$ and covariance $\bm{\Sigma}_{\mathbf{x}_0}$, stage costs $\ell$ and $\ell_f$, state dependent factor $\eta_t$, constraint sets $\mathbb{X}_t, \mathbb{U}_t$ $\forall t \in \{1,\ldots,T\}$, chance constraint probability $\epsilon$, ecdf confidence $\alpha$, tuning parameter $\delta$, decide if \textit{learning} should be carried out, the number of back-off iterations $n_b$ and the number of Monte Carlo simulations $S$ to estimate the back-offs.}
\item{Determine explicit back-off constraints using Algorithm 2.}
\item{Check final probabilistic value $\hat{\beta}_{lb}$ from Algorithm \ref{alg:back_off_algorithm} if it is close enough to $\epsilon$.}
\end{enumerate} 
\textit{Online Computations} \\ \For{$t=0,\ldots,T-1$}{
\begin{enumerate}
\item{Solve the MPC problem in Equation \ref{eq:nominalMPC} with the tightened constraint set from the \textit{Offline Computations}.}
\item{Apply the first control input of the optimal solution to the real plant.} 
\item{Measure the state $\mathbf{x}_t$ and update the GP plant model for learning GP NMPC.}
\end{enumerate}}
\end{algorithm2e}

Note the back-offs could be also updated online making use of the new initial conditions and updated prediction model in the case of "learning", which would lead to overall less conservativeness \citep{Paulson2018}. This would however require to carry-out the \textit{offline} calculations online, which is expensive, and not computationally desirable.

\section{Case study} \label{sec:case_study}
The case study utilized in this paper deals with the photo-production of phycocyanin synthesized by cyanobacterium \textit{Arthrospira platensis}. Phycocyanin is a high-value bioproduct and its biological function is to enhance the photosynthetic efficiency of cyanobacteria and red algae. It has been considered as a valuable compound because of its applications as a natural colorant to replace other toxic synthetic pigments in both food and cosmetic production. Furthermore, it has shown great promise for the pharmaceutical industry because of its unique antioxidant, neuroprotective, and anti-inflammatory properties. Using a simplified dynamic model we verify our GP NMPC algorithm by operating this process using a limited dataset. The GP NMPC problem is formulated with an economic objective aiming to directly maximize the bioproduct concentration of the final batch subject to two path constraints and one terminal constraint.   

\subsection{Semi-batch bioreactor model}
The simplified dynamic system consists of three ODEs describing the evolution of the concentration of biomass, nitrate, and bioproduct. The dynamic model is based on the Monod kinetics, which describes microorganism growth in nutrient sufficient cultures, where intracellular nutrient concentration is kept constant because of the rapid replenishment. We assume a fixed volume fed-batch. Control inputs are given by the light intensity ($I$) in $\upmu \text{mol.m}^{-2}$.s$^{-1}$ and nitrate inflow rate ($F_N$) in mg.L$^{-1}$.h$^{-1}$ \citep{DelRio-Chanona2015}. To capture the effects of light intensity on microalgae growth and bioproduction (photolimitation, photosaturation, and photoinhibition phenomena) the Aiba model is used \citep{Aiba1982}. The balance equations are given as follows:   

\begin{subequations} \label{eq:case_studyeq}
\begin{align}
    & \frac{dC_X}{dt} = u_m \cdot \frac{I}{I + k_s + \frac{I^2}{k_i}} \cdot C_X    \cdot \frac{C_N}{C_N + K_N} - u_d \cdot C_X, && {C_X}(0) = {C_X}_0 \\
    & \frac{dC_N}{dt} = -Y_{\frac{N}{X}} \cdot u_m \cdot \frac{I}{I + k_s + \frac{I^2}{k_i}} \cdot C_X \cdot \frac{C_N}{C_N + K_N} + F_N, && {C_N}(0) = {C_N}_0      \\
    & \frac{dC_{q_c}}{dt} = k_m \cdot \frac{I}{I + k_{sq} + \frac{I^2}{k_{iq}}} \cdot C_X - \frac{k_d C_{q_c}}{C_N + K_{Np}}, && {C_{q_c}}(0) = {C_{q_c}}_0 
\end{align}
\end{subequations}
where $C_X$ is the biomass concentration in $\text{g/L}$, $C_N$ is the nitrate concentration in $\text{mg/L}$, and $C_{q_c}$ is the phycocyanin (bioproduct) concentration in the culture in $\text{mg/L}$. The corresponding state vector and control vector are given by $\mathbf{x}=[C_X,C_N,C_{q_c}]^{\sf T}$ and $\mathbf{u}=[I,F_N]^{\sf T}$ respectively. The initial condition is denoted as $\mathbf{x}_0=[{C_X}_0,{C_N}_0,{C_{q_c}}_0]^{\sf T}$. The missing parameter values can be found in Table \ref{tab:par_casestudy}.

\begin{table}[H] 
\centering
\caption{Parameter values for ordinary differential equation system in Equation \ref{eq:case_studyeq}.}
\begin{tabular}{*{3}{l}} \label{tab:par_casestudy}
Parameter & Value & Units \\
\hline
$u_m$               & 0.0572  & $\text{h}^{-1}$                             \\
$u_d$               & 0.0     & $\text{h}^{-1}$                             \\
$K_N$               & 393.1   & mg.L$^{-1}$                                 \\
Y$_{\frac{N}{X}}$   & 504.5   & mg.g$^{-1}$                                 \\
$k_m$               & 0.00016 & mg.g$^{-1}$.h$^{-1}$                        \\
$k_d$               & 0.281   & h$^{-1}$                                    \\
$k_s$               & 178.9   & $\upmu \text{mol.m}^{-2}\text{.s}^{-1}$     \\
$k_i$               & 447.1   & $\upmu \text{mol.m}^{-2}\text{.s}^{-1}$     \\
$k_{sq}$            & 23.51   & $\upmu \text{mol.m}^{-2}\text{.s}^{-1}$     \\
$k_{iq}$            & 800.0   & $\upmu \text{mol.m}^{-2}\text{.s}^{-1}$     \\
$K_{Np}$            & 16.89   & mg.L$^{-1}$
\end{tabular}
\end{table}

\subsection{Problem set-up}
The time horizon $T$ was set to $12$ with an overall batch time of $240\text{h}$, and consequently the sampling time is $20 \text{h}$. Based on the dynamic system in Equation \ref{eq:case_studyeq} we define the objective and the constraints according to the general problem definition in Section \ref{sec:prob_def}. 
The measurement noise matrix $\bm{\Sigma}_{\bm{\upnu}}$ and disturbance noise matrix $\bm{\Sigma}_{\bm{\upomega}}$ were set to:
\begin{equation}
    \bm{\Sigma}_{\bm{\upnu}} = \diag(4 \times 10^{-4},0.1,1 \times 10^{-8}),  \quad \bm{\Sigma}_{\bm{\upomega}} = \diag(4 \times 10^{-4},0.1,1 \times 10^{-8})  
\end{equation}

The mean $\bm{\upmu}_{\mathbf{x}_0}$ and covariance $\bm{\Sigma}_{\mathbf{x}_0}$  of the initial condition are given by:
\begin{equation}
    \bm{\upmu}_{\mathbf{x}_0} = [1.,150,0.]^{\sf T}, \quad  \bm{\Sigma}_{\mathbf{x}_0} = \diag(1 \times 10^{-3},22.5,0)
\end{equation}

The control algorithm aims to maximize the amount of bioproduct produced $C_{q_c}$ with a penalty on the change of control actions. The corresponding stage and terminal cost can be stated as follows: 
\begin{subequations}
\begin{align}
     & \ell(\mathbf{x}_t,\mathbf{u}_t) = \bm{\Delta}_{\mathbf{u}_t}^{\sf T} \mathbf{R} {\bm{\Delta}}_{\mathbf{u}_t} \\
     & \ell_f(\mathbf{x}_T) = -{C_{q_c}}_T
\end{align}
\end{subequations}
where ${\bm{\Delta}}_{\mathbf{u}_t} = \mathbf{u}_{t} - \mathbf{u}_{t-1}$ and $\mathbf{R} = \diag(3.125 \times 10^{-8},3.125 \times 10^{-6})$. The overall objective is then defined by Equation \ref{eq:objdef}. 

For the case of state dependency we set all $\eta_i$ to zero except $\eta_0$, see Equation $\ref{eq:nominalMPC}$. The value of $\eta_0$ was set to 15. Note that for these factors to work properly it is important to normalize the data as we did in this case study. 

There are two path constraints in the problem. The amount of nitrate is constrained to remain below $800$ mg/L, while the ratio of bioproduct to biomass may not exceed $11.0$ mg/g for high density biomass cultivation. These constraints can be stated as:
\begin{subequations}
\begin{align}
    & g_1^{(t)}(\mathbf{x}_t) = {C_N}_t - 800 \leq 0 && \forall t \in \{0,\ldots,T\}  \\
    & g_2^{(t)}(\mathbf{x}_t) = {C_{q_c}}_t - 0.011 {C_X}_t \leq 0 && \forall t \in \{0,\ldots,T\}
\end{align}
\end{subequations}

Lastly, there is a terminal constraint on nitrate to reach a final concentration of below $150$ mg/L:
\begin{align}
    & g_3^{(T)}(\mathbf{x}_T) =  {C_N}_T - 150 \leq 0, \quad g_3^{(t)}(\mathbf{x}_T) = 0 \; \forall t \in \{0,\ldots,T-1\} 
\end{align}

The maximum probability for violating the joint chance constraint was set to $\epsilon=0.1$. The control inputs light intensity and nitrate inflow rate are constrained as follows:   

\begin{subequations}
\begin{align}
    & 120 \leq I_t \leq 400 && \forall t \in \{0,\ldots,T\} \\
    & 0 \leq {F_N}_t \leq 40    && \forall t \in \{0,\ldots,T\} 
\end{align}
\end{subequations}

For the back-off iterations we employed $S=1000$ MC iterations with the initial back-offs computed according to Equation \ref{eq:update_rule} with $\delta=0.1$ and $\alpha=0.01$. The maximum number of back-off iterations for the bisection algorithm was set to $n_b=16$.   

\subsection{Implementation details and initial dataset generation}
The optimization problem for the GP NMPC in Equation \ref{eq:nominalMPC} is solved using Casadi \citep{Andersson2019} to obtain the gradients of the problem using automatic differentiation in conjunction with IPOPT \citep{Wachter2006}. The "real" plant model was simulated using IDAS \citep{Hindmarsh2005}. In the next section, different variations of the proposed algorithm are presented, for which two different type of datasets were collected. For the first type of dataset we designed the entire input data matrix $\mathbf{Z}$ according to a Sobol sequence \citep{Sobol2001} in the range $\mathbf{z}_i \in [0,20] \times [50,800] \times [0,0.18] \times [120,400] \times [0,40]$. The ranges were chosen for the data to cover the expected operating region. The corresponding outputs $\mathbf{Y}$ were then obtained from the IDAS simulation of the system perturbed by Gaussian noise as defined in the problem setup. In the second approach only the control inputs were set according to the Sobol sequence in the range $\mathbf{u}_i \in [120,400] \times [0,40]$ and the corresponding states $\mathbf{Y}$ were obtained from the trajectories of the "real" system perturbed by noise using samples of the initial condition and the time horizon as defined in the problem setup based on these control inputs. The system was simulated in "open-loop" using these control actions, i.e. without any feedback controller present. For both datasets the input data $\mathbf{Z}$ and output data $\mathbf{Y}$ are normalized to zero mean and a standard deviation of one. The reason we use two different types of datasets is to highlight the advantages of accounting for state dependency in two of the algorithm variations. In the first dataset the data is relatively evenly distributed and hence considering the state dependency of the uncertainty to avoid regions with high data sparsity has essentially no effect, while in the second approach there are clearly defined trajectories that can be followed by accounting for state dependency.  

\section{Results and discussions} \label{sec:results_discussions}
In this section we present and discuss the results from the case study described in the previous section. For comparison purposes we compare six different variations of the proposed GP NMPC approach, which are as follows:
\begin{itemize}
    \item{GP NMPC 50, 60, 100: GP NMPC approach without learning and without taking into account state dependency for dataset sizes of 50, 60, and 100 points using the first type of dataset.}
    \item{GP NMPC learning 50: GP NMPC approach with learning and without state dependency for a dataset size of 50 points, which will be compared to the above case of 50 data points without learning. The first type dataset is utilized.}
    \item{GP NMPC SD/NSD 50: GP NMPC approach with and without accounting for the state dependency for a dataset size of 50 points employing the second type of dataset.}
\end{itemize}

\begin{figure}[H] \centering
   \includegraphics[width=0.75\textwidth]{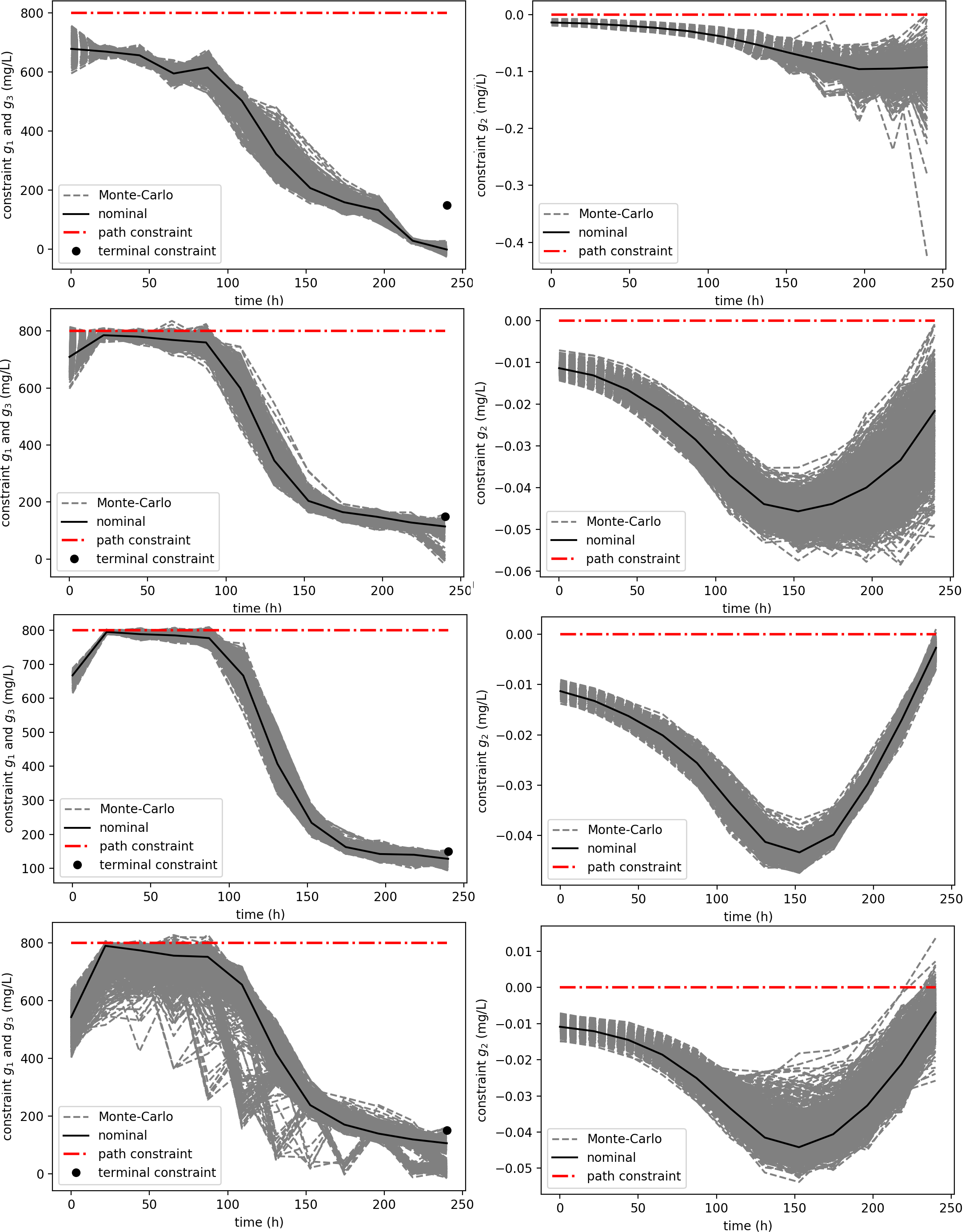} \caption{The 1000 MC trajectories at the final back-off iteration of the nitrate concentration constraints $g_1$ and $g_3$ (LHS) and the ratio of bioproduct to biomass constraint $g_2$ (RHS) for from top to bottom GP NMPC 50, 60, 100 and GP NMPC 50 learning.}
  \label{fig:MC_con_50_60_100_learning}
\end{figure}

\begin{figure}[H] \centering
   \includegraphics[width=0.75\textwidth]{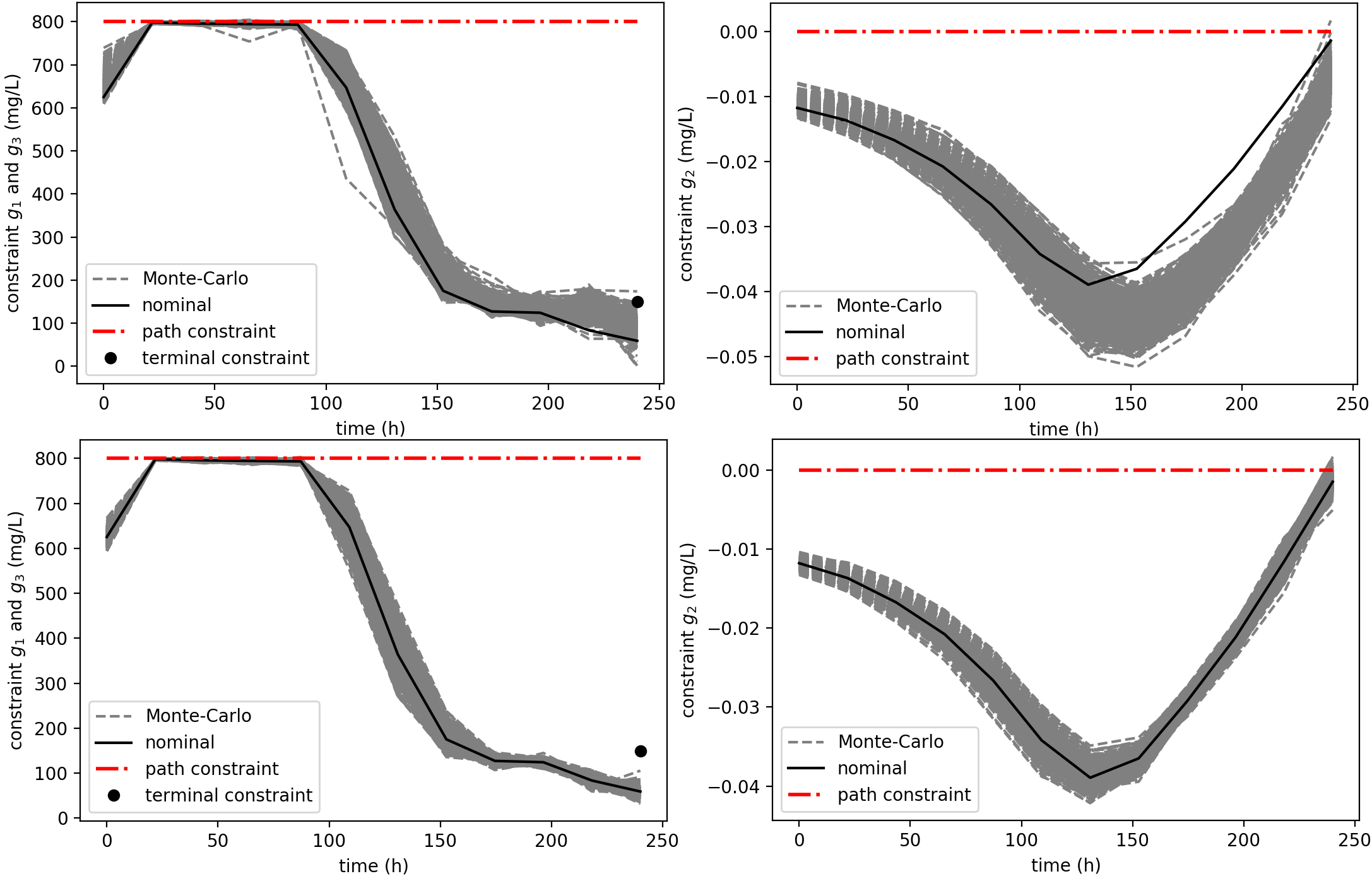} \caption{The 1000 MC trajectories at the final back-off iteration of the nitrate concentration constraints $g_1$ and $g_3$ (LHS) and the ratio of bioproduct to biomass constraint $g_2$ (RHS) for GP NMPC NSD 50 (top) and GP NMPC SD 50 (bottom).}
  \label{fig:MC_con_SD_NSD}
\end{figure}

\begin{figure}[H]
\centering
   \includegraphics[width=0.8\textwidth]{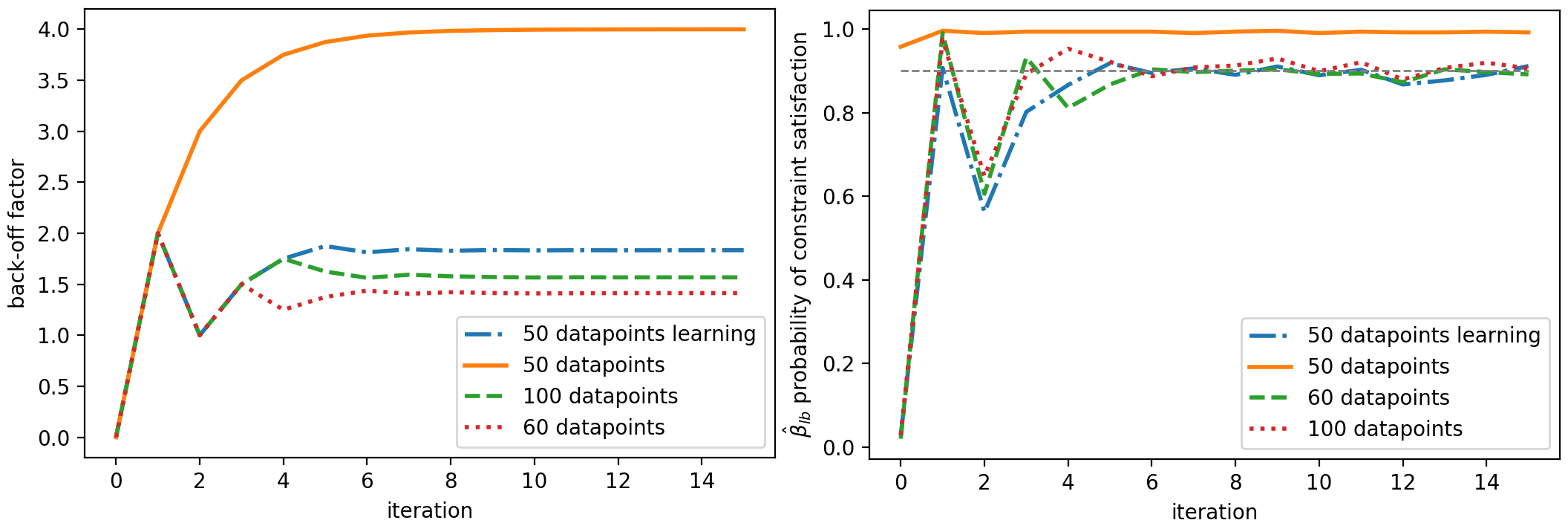} \caption{Plots of evolution of the back-off factor and the probability of constraint satisfaction $\hat{\beta}_{lb}$ over the 16 back-off iterations for GP NMPC 50, 60, 100, and GP NMPC learning 50.}
  \label{fig:back_off_beta_50_60_100_learning}
\end{figure}

\begin{figure}[H] \centering
   \includegraphics[width=0.8\textwidth]{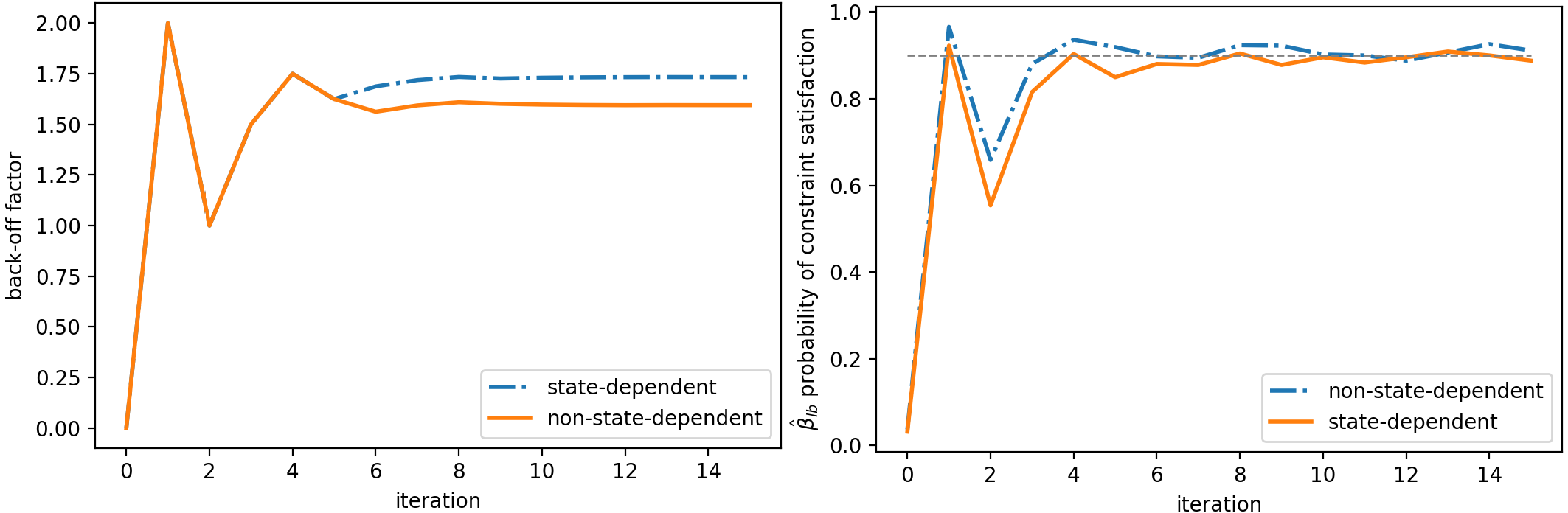} \caption{Plots of evolution of the back-off factor and the probability of constraint satisfaction $\hat{\beta}_{lb}$ over the 16 back-off iterations for GP NMPC SD 50 and GP NMPC NSD 50.}
  \label{fig:back_off_beta_SD_NSD}
\end{figure}

\begin{figure}[H] \centering
   \includegraphics[width=0.8\textwidth]{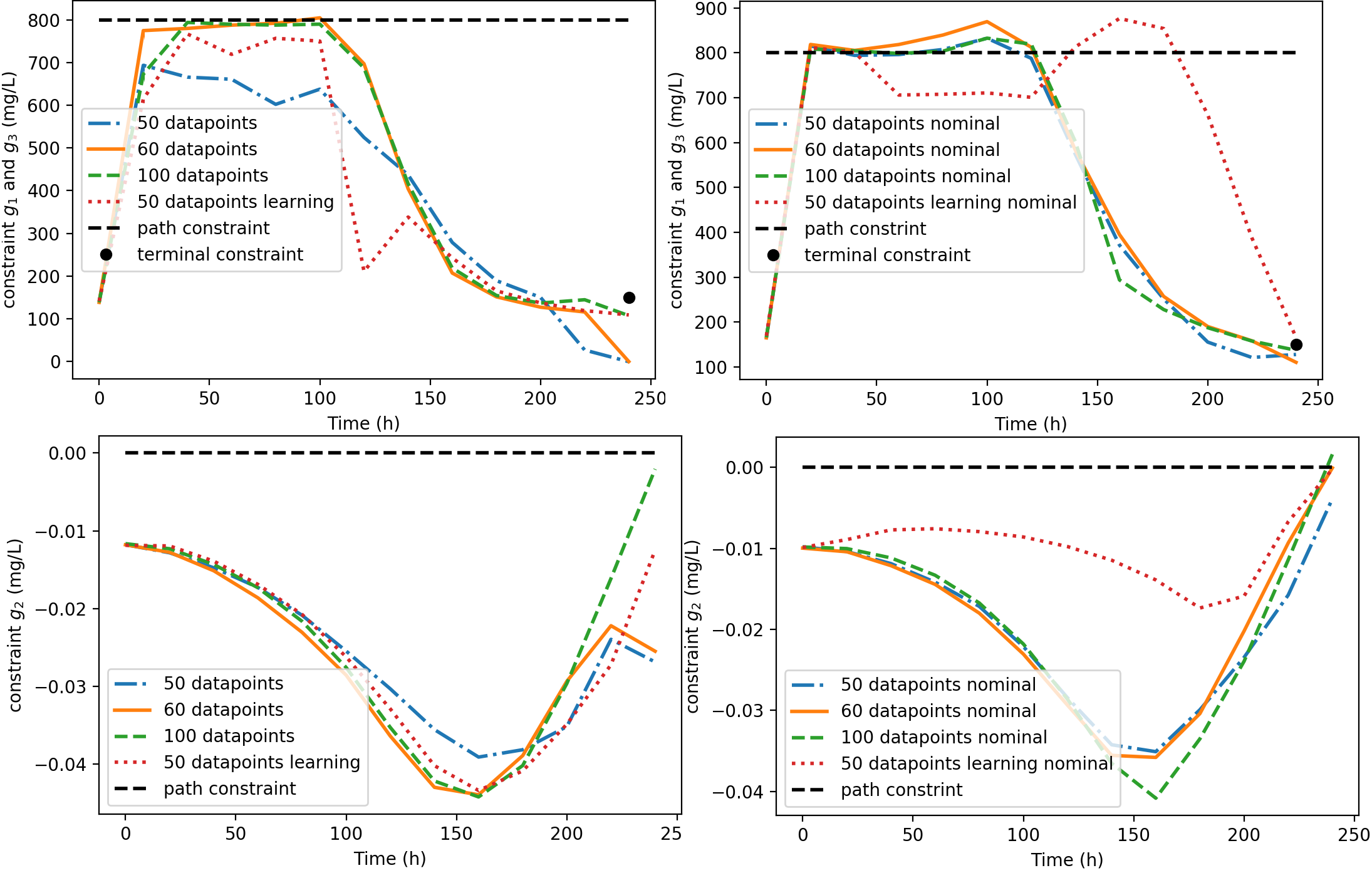} \caption{Trajectories of the nitrate concentration constraints $g_1$ and $g_3$ (top) and the ratio of bioproduct to biomass constraint $g_2$ (bottom) for the GP NMPC 50, 60, 100, and GP NMPC 50 learning applied to the "real" plant model with the final tightened constraint set on the LHS and with no back-off constraints on the RHS referred to as $\textit{nominal}$.}
  \label{fig:plant_50_60_100_learning}
\end{figure}

\begin{figure}[H] \centering
   \includegraphics[width=0.8\textwidth]{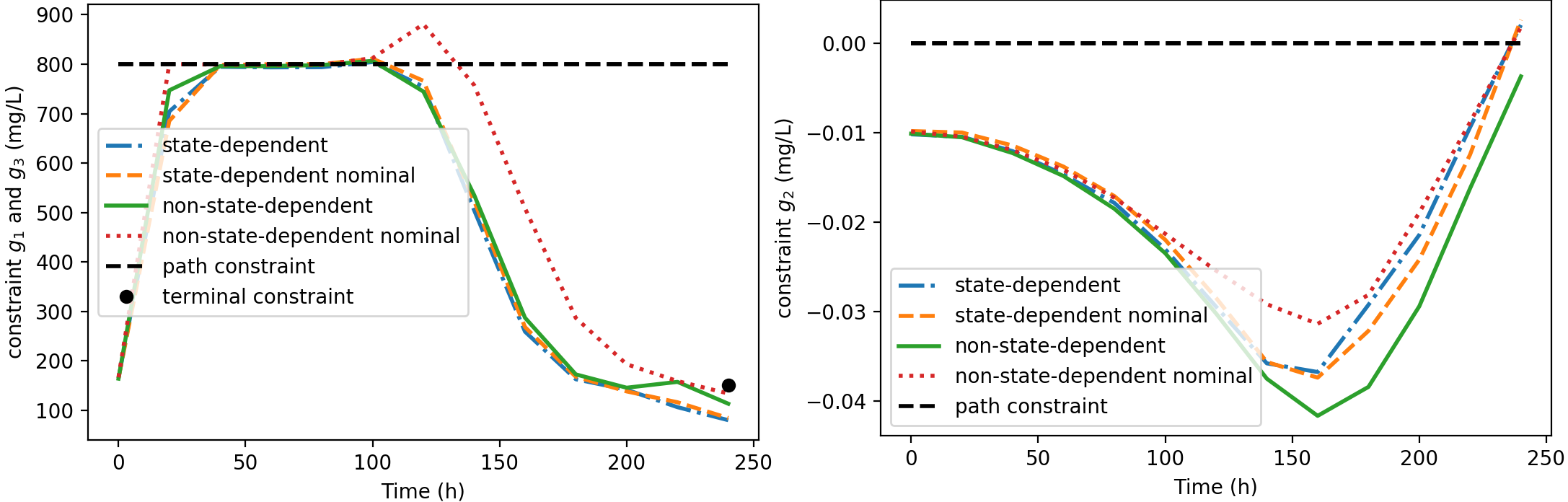} \caption{Trajectories of the nitrate concentration constraints $g_1$ and $g_3$ (LHS) and the ratio of bioproduct to biomass constraint $g_2$ (RHS) for the GP NMPC 50, 60, 100, and GP NMPC 50 learning applied to the "real" plant model with the final tightened constraint set and with no back-off constraints referred to as $\textit{nominal}$.}
  \label{fig:plant_SD_NSD}
\end{figure}

\begin{figure}[H] \centering
   \includegraphics[width=0.8\textwidth]{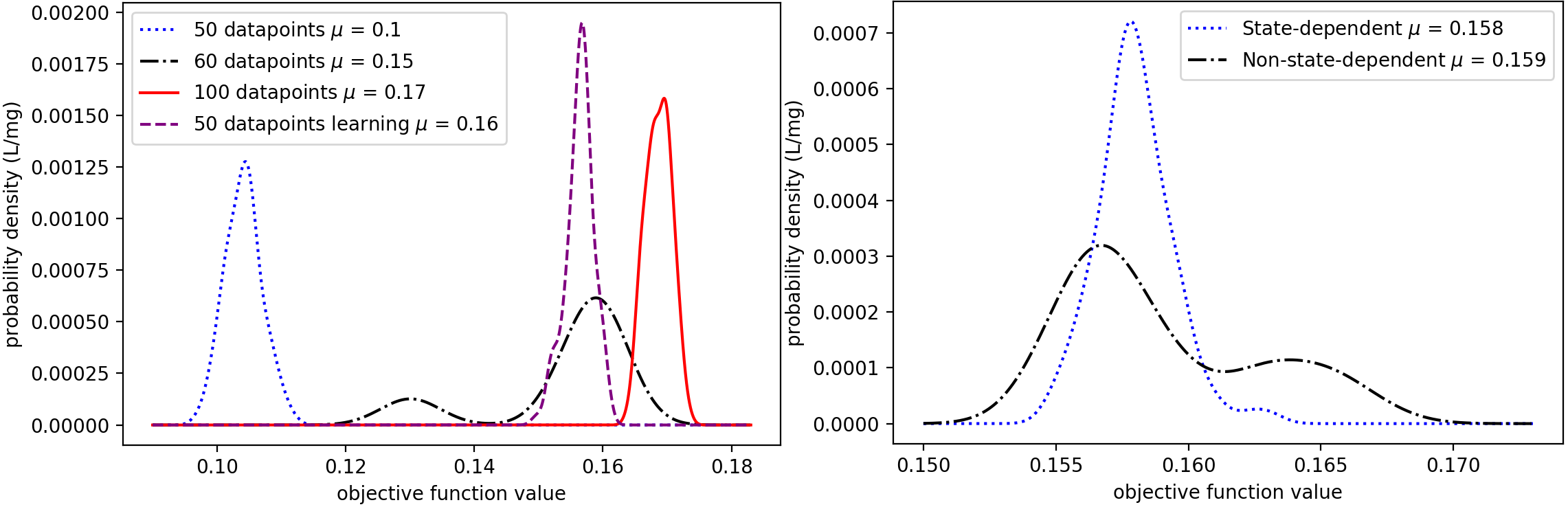} \caption{Probability density function for the "real" plant objective values for all variations of the GP NMPC algorithm.}
  \label{fig:objective_all}
\end{figure}

\begin{figure}[H] \centering
   \includegraphics[width=0.8\textwidth]{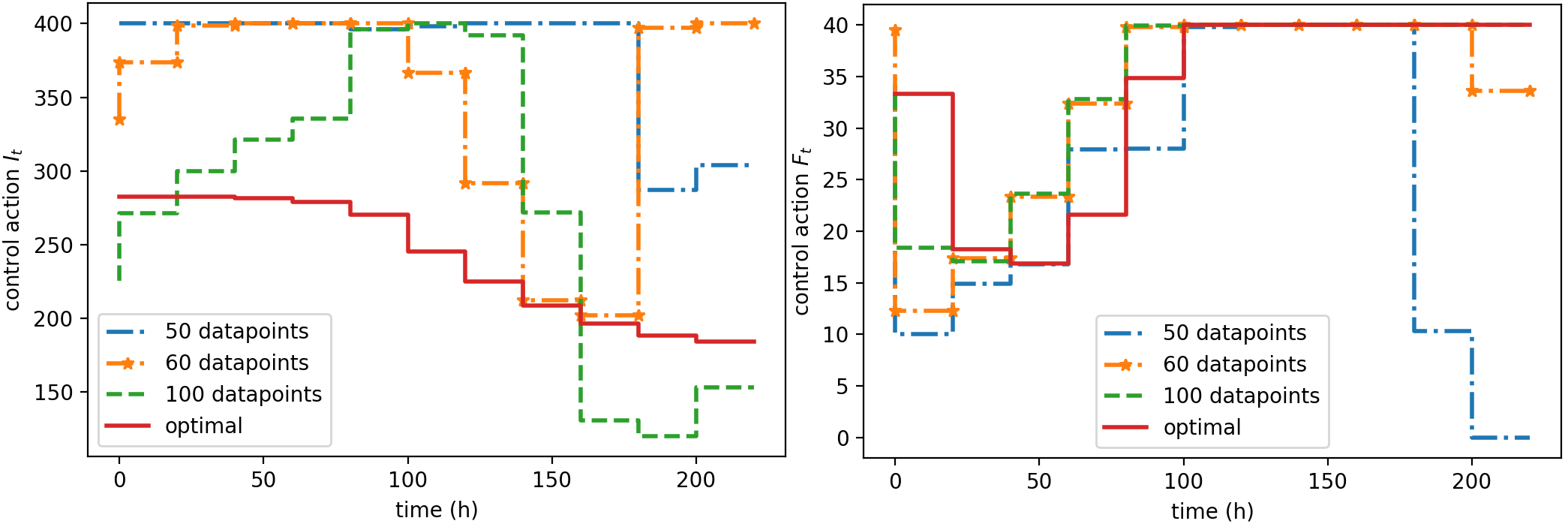} \caption{Example control trajectories for the light intensity on the LHS and the nitrate flow rate on the RHS based on GP NMPC 50, 60, 100. The red line represents the optimal control trajectories ignoring the noise present in the process. }
  \label{fig:control_trajectories_50_60_100}
\end{figure}

\begin{figure}[H] \centering
   \includegraphics[width=0.8\textwidth]{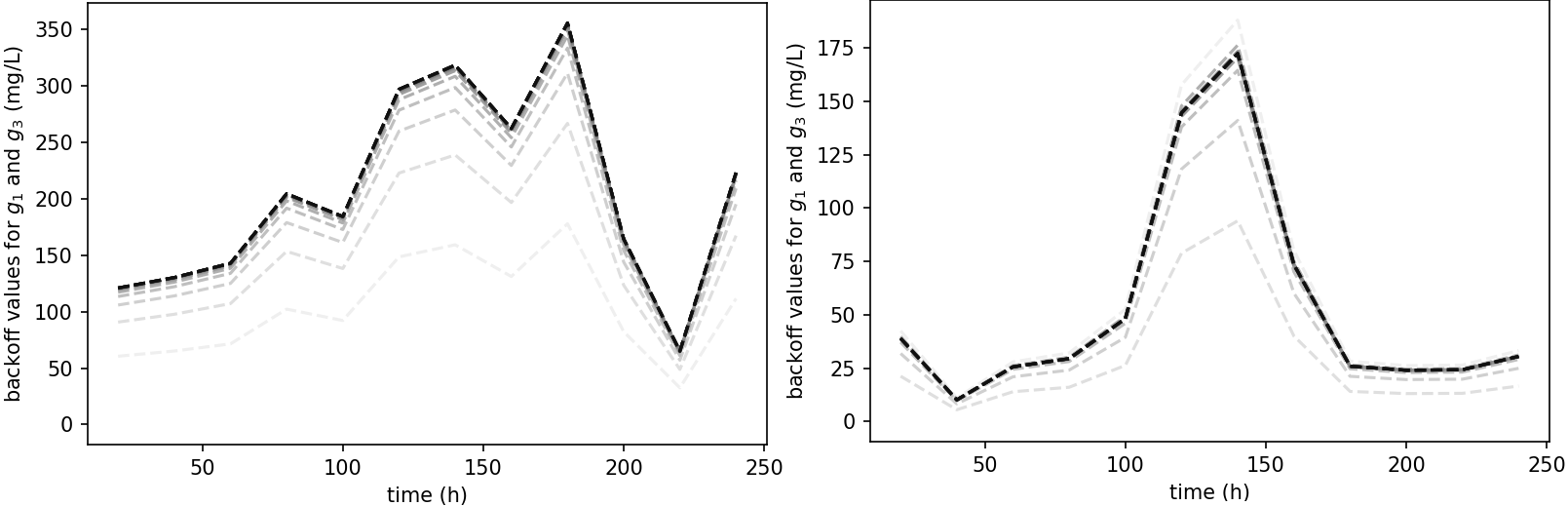} \caption{Example back-off values for the nitrate concentration constraints $g_1$ and $g_3$ of GP NMPC 50 (LHS) and GP NMPC 50 learning (RHS). The lines are plotted over the 16 back-off iterations, which are faded out towards earlier iterations.}
  \label{fig:example_back_off_learning_nonlearning}
\end{figure}

In addition, we compare the approaches to a nominal NMPC algorithm based on the GP model to show the importance of employing back-offs to prevent constraint violations, i.e. we run the GP NMPC on the "real" plant model, while setting the back-offs to zero. The results of the outlined runs are summarized in Figures \ref{fig:back_off_beta_50_60_100_learning}-\ref{fig:example_back_off_learning_nonlearning}, and in Tables \ref{tab:mean_back_off}-\ref{tab:probability_OCP_time}. In Figures \ref{fig:back_off_beta_50_60_100_learning}-\ref{fig:back_off_beta_SD_NSD} we show the evolution of the back-off factor and the probability of constraint satisfaction $\hat{\beta}_{lb}$ over the 16 back-off iterations from Algorithm $\ref{alg:back_off_algorithm}$. The next two Figures \ref{fig:MC_con_50_60_100_learning}-\ref{fig:MC_con_SD_NSD} show the 1000 MC trajectories of the constraints with a line to highlight the nominal prediction of the GP NMPC. Next the GP NMPC was applied to the "real" plant with back-offs from the final iteration and without back-offs referred to as $\textit{nominal}$ as shown in Figures \ref{fig:plant_50_60_100_learning}-\ref{fig:plant_SD_NSD}. Figure \ref{fig:objective_all} shows the probability density function of the objective values obtained from the "real" plant, where in the figure larger objective values correspond to better objective values. Figure \ref{fig:control_trajectories_50_60_100} shows representative control trajectories for GP NMPC 50, 60, 100 compared with the optimal trajectory obtained from solving the OCP of the "real" plant ignoring uncertainties. Figure \ref{fig:example_back_off_learning_nonlearning} shows the back-off values for the nitrate constraints $g_1$ and $g_2$ for GP NMPC 50 and GP NMPC 50 learning. Lastly, Table $\ref{tab:mean_back_off}$ shows the mean values for the back-offs averaged over time for the final back-off iteration, while Table $\ref{tab:probability_OCP_time}$ shows the attained probability of satisfaction $\hat{\beta}_{lb}$ together with the average computational times for solving a single GP NMPC optimization problem. We can draw the following conclusions from these results:

\begin{itemize}
\item{Figures \ref{fig:back_off_beta_50_60_100_learning}-\ref{fig:back_off_beta_SD_NSD} and Table $\ref{tab:mean_back_off}$ show that apart from GP NMPC 50 the other variations reach the required $\hat{\beta}_{lb}$ and hence successfully converge to a reasonable back-off factor. For these, as expected, a low back-off value leads to too low $\hat{\beta}_{lb}$ values near zero, while too high back-off values lead to too high $\hat{\beta}_{lb}$ values. The value of $\hat{\beta}_{lb}$ does vary by $\pm 0.01$ even on convergence, which is due to the randomness of the MC samples. Nonetheless, since  $\hat{\beta}_{lb}$ is a sample robust value, it is high enough if at least $0.9$ is reached once over the 16 iterations, which is the case for all of them. Note that GP NMPC 50 does not converge, since even without back-offs the NMPC remains feasible for all MC trajectories and hence the bisection procedure fails. The performance of GP NMPC 50 is however also by far the worst. This is due to insufficient amounts of data in crucial areas for the control problem.}

\item{From GP NMPC 50, 60 to 100 the objective values steadily increase and hence improve with increased number of data points as shown in Figure \ref{fig:objective_all}. This is as expected, since more data points should lead to a more accurate GP plant model and hence more optimal control actions. Lastly, a more accurate GP plant model should also require less conservative back-offs. For the constraint $g_2$ the mean of the back-off values steadily decreases from $0.022$ mg/L for GP NMPC 50 to $0.003$ for GP NMPC 100 as shown in Table \ref{tab:mean_back_off}. For the constraints $g_1/g_3$ on the other hand the mean of the back-offs decreases dramatically from GP NMPC 50 with a value of $250$mg/L to a value of $34.2$mg/L for GP NMPC 60, while slightly increasing again for GP NMPC 100 to $38.8$ mg/L. This is further illustrated in Figure \ref{fig:MC_con_50_60_100_learning}, for which the spread of the trajectories decreases steadily from GP NMPC 50, 60 to 100. All in all, larger datasets lead to improved solutions.}

\item{The learning approach GP NMPC 50 learning leads to a reasonable solution with an objective value that is on average higher and therefore an improvement over GP NMPC 60 as can be seen in Figure \ref{fig:objective_all}. Further, the sharper peak of the objective value suggests a more reliable performance. In contrast, GP NMPC 50 without learning is unable to determine a good solution and therefore has an objective value that is considerably worse than the remaining scenarios with an objective value that is on average over 30$\%$ lower. GP NMPC 50 is also unable to reach a $\hat{\beta}_{lb}$ value of $0.9$ and has instead a much higher value as shown in Table \ref{tab:probability_OCP_time}, but at the expense of performance. This is highlighted in Figure \ref{fig:control_trajectories_50_60_100} in which the control trajectory of GP NMPC 50 decreases the nitrate flowrate early to satisfy the terminal constraint, which leads to a sub optimal solution. This is believed to be due to the high uncertainty of the $g_1/g_2/g_3$ constraint trajectories of GP NMPC 50, which can be seen by the large spread of the trajectories in Figure \ref{fig:MC_con_50_60_100_learning}. GP NMPC 50 learning on the other hand has a spread of the constraints $g_1/g_2/g_3$ that is significantly less than GP NMPC 50. This is further highlighted by the considerably higher back-off values of GP NMPC 50 compared to GP NMPC 50 learning as can be seen in Table \ref{tab:mean_back_off}, where the $g_1/g_3$ back-off values are nearly 400$\%$ larger, while the $g_2$ back-off values are over $300\%$ larger. In conclusion, accounting for online learning has lead to a significantly better solution, although it should be noted that at larger datasets the effect is nearly negligible. }

\item{GP NMPC 50 can be seen to be less erratic than GP NMPC 50 learning in Figure \ref{fig:MC_con_50_60_100_learning}, which is however expected. GP NMPC 50 uses the same prediction model throughout and hence the control inputs are only influenced by changes in the initial condition. For GP NMPC 50 learning on the other hand the prediction model changes at each sampling point and hence this leads to more irregular behaviour, which can be seen by the increased oscillations. This then leads to overall improved control actions and performance due to exploiting the new data available, as can be seen in Figure \ref{fig:objective_all}.}

\item{It can be seen in Figure \ref{fig:MC_con_SD_NSD} that GP NMPC NSD 50 has a larger spread of trajectories than GP NMPC SD 50. This is as expected, since GP NMPC SD 50 aims directly in its objective to minimize uncertainty. Consequently, the mean of the back-off values of constraints $g_1/g_3$ and $g_2$ in Table $\ref{tab:mean_back_off}$ are over $50\%$ larger and over $100\%$ greater than those for GP NMPC SD 50, respectively. Nonetheless, the attained average objective of GP NMPC SD 50 is marginally lower and hence worse as can be seen in Figure \ref{fig:objective_all}. This is expected, since accounting for state dependency reduces conservativeness, but may also lead to a sub optimal solution due to conflicting objectives.}

\item{In Table \ref{tab:probability_OCP_time} the average computational times of a single GP NMPC evaluation are shown, which range from $135$ms to $48$ms. Overall, it can be seen that by far the largest computational times are attributed to GP NMPC 100, which is reasonable since the complexity of GP plant models grows exponentially with the number of data points. GP NMPC 50 and GP NMPC 50 learning can be seen however to have higher computational times, which is due to a more complex optimization problem from the reduced amount of data. This is further highlighted by GP NMPC 50 SD and GP NMPC 50 NSD attaining the lowest average computational times due to the second type of dataset leading to easier optimization solutions. In Table \ref{tab:probability_OCP_time} we can further see that the computational times required for a single back-off iteration is for all variations nearly solely determined by the GP NMPC evaluation time. In addition, in Table \ref{tab:probability_OCP_time} we show the ecdf value $\hat{\beta}$ for the final back-off iteration. It can be seen that these probabilities are substantially higher than the required probability of $0.9$ ranging from $0.91$ to $0.93$ for the converged solutions, which is due to the conservativeness of the probabilistic lower bound used. This leads to higher back-off values than required and therefore worse objective values. This conservativeness can be reduced by setting $\alpha$ to a higher value.}

\item{Lastly, in Figures \ref{fig:plant_50_60_100_learning}-\ref{fig:plant_SD_NSD} trajectories of the constraints are shown by applying the GP NMPC variants to the "real" plant. It can be seen that the $\textit{nominal}$ variations of GP NMPC 50, 60, 100, and GP NMPC 50 learning with back-offs set to zero violate the nitrate constraint $g_1$ to remain below $800$mg/L by a substantial amount of up to $50$mg/L for GP NMPC 50 learning and GP NMPC 60. With back-offs on the other hand the approaches remain feasible throughout the run, which illustrates the importance of employing back-offs. GP NMPC NSD 50 $\textit{nominal}$ can be also seen to violate the nitrate constraint $g_1$ by $50$mg/L, while GP NMPC SD 50 $\textit{nominal}$ does not violate this constraint. This is likely due to GP NMPC SD 50 $\textit{nominal}$ following a feasible trajectory in the dataset. GP NMPC NSD 50 with back-offs remains feasible. Overall, it can be seen that back-offs are important to achieve feasibility given the presence of plant-model mismatch.}

\end{itemize}

\begin{table}[H]
\centering
\caption{The mean of the back-off values for the nitrate concentration constraints $g_1/g_3$ and the ratio of bioproduct to biomass constraint $g_2$ from the final back-off iteration.}
\begin{tabular}{*{3}{l}}
Algorithm variation & Mean back-off $g_1/g_3$ (mg/L) & Mean back-off $g_2$ (mg/L) \\
\hline
GP NMPC 50          & 205.8 & 0.022    \\
GP NMPC 60          & 34.2  & 0.008    \\
GP NMPC 100         & 38.8  & 0.003    \\
GP NMPC learning 50 & 54.0  & 0.007    \\
GP NMPC 50 SD       & 23.3  & 0.002    \\
GP NMPC 50 NSD      & 36.2  & 0.004
\end{tabular}
\label{tab:mean_back_off}
\end{table}

\begin{table}[H]
\centering
\caption{Lower bound on the probability of satisfying the joint constraint $\hat{\beta}_{lb}$, the approximate satisfaction probability from the final simulation $\hat{\beta}$, average computational times to solve a single optimal control problem (OCP) for the GP NMPC, and the average computational time required to complete one back-off iteration.}
\begin{tabular}{*{5}{l}}
Algorithm variation & Probability $\hat{\beta}_{lb}$ & Probability $\hat{\beta}$ & OCP time (ms) & Back-off iteration time (s) \\
\hline
GP NMPC 50          & 0.99  & 1.00 & 65  & 782  \\
GP NMPC 60          & 0.89  & 0.91 & 54  & 753  \\
GP NMPC 100         & 0.91  & 0.93 & 135 & 1626 \\
GP NMPC learning 50 & 0.91  & 0.93 & 69  & 825  \\
GP NMPC 50 SD       & 0.89  & 0.91 & 48  & 574  \\
GP NMPC 50 NSD      & 0.91  & 0.93 & 49  & 584
\end{tabular}
\label{tab:probability_OCP_time}
\end{table}

\section{Conclusions} \label{sec:conclusions}
In conclusion, a new approach is proposed for finite-horizon control problems using NMPC in conjunction with GP state space models. The method utilizes the probabilistic nature of GPs to sample deterministic functions of possible plant models. Tightened constraints using explicit back-offs are then determined, such that the closed-loop simulations of these possible plant models are feasible to a high probability. In addition, it is shown how probabilistic guarantees can be derived based on the number of constraint violations from the simulations. Furthermore, it is shown that online learning and state dependency of the uncertainty can be taken into account explicitly in this method, which leads to overall less conservativeness. Moreover, the computational times are shown to be relatively low, since constraint tightening is performed offline. Finally, through the comprehensive semi-batch bioprocess case study, the efficiency and potential of this method for the optimisation of complex stochastic systems (e.g. biological processes) is well demonstrated.

\section*{Acknowledgements}
This project has received funding from the European Union’s Horizon 2020 research and innovation programme under the Marie Sklodowska-Curie grant agreement No 675215.

\section*{References}
\bibliography{Backoff_GP_JPC}
\bibliographystyle{model4-names}\biboptions{authoryear}

\appendix

\section{Recursive inverse matrix update} \label{app:recursive_update}
In this paper we are often concerned with inverting matrices as shown in Section \ref{sec:Gaussian_processes_online_learning} in Equation \ref{eq:update_cov}:
\begin{align}
        & \bm{\Sigma}^{+-1}_{\mathbf{Y}} = \begin{bmatrix}
    \bm{\Sigma}_{\mathbf{Y}} & \mathbf{k}^{\sf T}(\mathbf{z}^+)  \\
    \mathbf{k}(\mathbf{z}^+) & k(\mathbf{z}^+,\mathbf{z}^+) (+ \sigma_{\nu}^2)
\end{bmatrix}^{-1}
\end{align}

This is an recursive update formula and hence $\bm{\Sigma}_{\mathbf{Y}}^{-1}$ is known \textit{a priori}. In this Section we introduce a recursive formula to exploit this fact to obtain $\bm{\Sigma}^{+-1}_{\mathbf{Y}}$ in a cheaper fashion taken from \citet{Strassen1969} adjusted to our case. The following quantities need to be computed to obtain $\bm{\Sigma}^{+-1}_{\mathbf{Y}}$:  

\begin{subequations}
\begin{align}
    & \RN{1}     && = \mathbf{k}^{\sf T}(\mathbf{z}^+) \bm{\Sigma}_{\mathbf{Y}}^{-1}                   \\
    & \RN{2}     && = \mathbf{k}^{\sf T}(\mathbf{z}^+) \RN{1}^{\sf T}                                  \\
    & \mathbf{C}_{12}     && = \RN{1}^{\sf T} \times \RN{2}                                             \\
    & C_{11}     && = \bm{\Sigma}_{\mathbf{Y}}^{-1} - \RN{1}^{\sf T} \times \mathbf{C}_{12}^{\sf T}  \\
    & C_{22}     && = -\left(\RN{2} - k(\mathbf{z}^+,\mathbf{z}^+) (+ \sigma_{\nu}^2) \right)^{-1}
\end{align}
\end{subequations}

The inverted matrix is then given by:
\begin{equation}
    \bm{\Sigma}^{+-1}_{\mathbf{Y}} = \begin{bmatrix}
    C_{11} & \mathbf{C}_{12}  \\
    \mathbf{C}_{12}^{\sf T} & C_{22}
\end{bmatrix}
\end{equation}
\end{document}